\documentclass{amsart}
\usepackage{amsthm}
\usepackage{geometry}
\usepackage{amsmath,amssymb}
\usepackage{graphicx}
\newtheorem{theorem}{Theorem}[section]
\newtheorem*{Thm1prime}{Theorem \ref{main1}'}
\newtheorem{lemma}[theorem]{Lemma}
\newtheorem{proposition}[theorem]{Proposition}

\newtheorem{definition}[theorem]{Definition}

\newtheorem{remark}[theorem]{Remark}

\newcommand{\N}{{\mathbb N}}
\newcommand{\R}{{\mathbb R}}

\newcommand{\Z}{{\mathbb Z}}

\begin{document}

\title{The existence of small prime gaps in subsets of the integers}

\author{Jacques Benatar}
\address{UCLA Department of Mathematics, Los Angeles, CA 90095-1555.}
\email{jbenatar@math.ucla.edu}
 \maketitle
 
\begin{abstract}
We consider the problem of finding small prime gaps in various sets $\mathcal{C} \subset \N$. Following the work of Goldston-Pintz-Y{\i}ld{\i}r{\i}m, we will consider collections of natural numbers that are well-controlled in arithmetic progressions. Letting $q_n$ denote the $n$-th prime in $\mathcal{C}$, we will establish that for any small constant $\epsilon>0$, the set $\left\{ q_n| q_{n+1}-q_n \leq \epsilon \log n \right\}$ constitutes a positive proportion of all prime numbers. Using the techniques developed by Maynard and Tao we will also demonstrate that $\mathcal{C}$ has bounded prime gaps. Specific examples, such as the case where $\mathcal{C}$ is an arithmetic progression have already been studied and so the purpose of this paper is to present results for general classes of sets.
\end{abstract}

\section{Introduction and framework}
One of the most famous unsolved problems in Number Theory is the so-called Twin Prime Conjecture, which posits the existence of infinitely many pairs of primes $(p,p')$ for which $p-p'=2$. Throughout the past century a great amount of work has been done with regards to this conjecture and we refer the reader to \cite{3} and \cite{Sound} for some historical background on the subject. In this article we will build upon the methods developed by Goldston, Pintz, Y{\i}ld{\i}r{\i}m and more recently Maynard and Tao. Letting $p_n$ denote the $n$-th prime number, it was shown in \cite{May} that 
$$\liminf_{n \rightarrow \infty} (p_{n+m} - p_n) \ll e^{4m+ \varepsilon} $$
for any $\varepsilon>0$ and, regarding the frequency of small prime gaps, it was demonstrated in \cite{PTIV} that 
\[ |\left\{ p_n \leq N | \  p_{n+1}-p_n \leq \eta \log N  \right\}| \geq C \frac{N}{\log N}\]
for some constant $C$ depending on $\eta$ and any natural number $N$.
In this paper we will investigate which conditions ensure that a set $\mathcal{C} \subset \N$ will also have many small prime gaps.
\begin{definition}
Let $q_n$ denote the $n$-th prime number in $\mathcal{C}$ and for $\eta>0 $, write
\[ \pi_\mathcal{C}(N;\eta):=|\left\{ q_n \leq N | \  q_{n+1}-q_n \leq \eta \log N  \right\}|. \]
We say that $\mathcal{C}$ frequently contains $\eta$-small prime gaps if there exists a constant $C$, depending on $\eta$, such that
\[ \pi_\mathcal{C}(N;\eta) \geq C \frac{N}{\log N} \]
for all positive integers $N\geq 1$. We also define the quantity
$$H_{\mathcal{C}}(m):= \liminf_{n \rightarrow \infty} (q_{n+m} - q_n).$$
\end{definition}
From a probabilistic point of view, let us first show that a random subset of the primes, with positive relative density $\rho$, frequently contains $\eta$-small prime gaps (for any $\eta$). To avoid issues of independence, we will work with the sets $\mathcal{P}_\mathcal{C}(N;\eta)^{(i)}:=\left\{ p_n \leq N | p_n,p_{n+1} \in \mathcal{C}, n \equiv i \bmod 2 \text{ and }  p_{n+1}-p_n \leq \eta \log N  \right\}$ for $i=1,2$. For fixed $N$, denote $\pi_\mathcal{C}^{*}(N;\eta)$ for the cardinality of the largest of these two sets. 
\begin{lemma}
Let $B$ be a random subset of the primes defined by letting the independent events $p \in B$ occur with probability $\mathbb{P}(p \in B)=\rho$. Writing $\lambda(N,\eta):=\rho^2 \pi_{\N}^{*}(N;\eta) \geq \rho^2 \pi_{\N}(N;\eta)/2$ one has that
\begin{align}\label{heuristics}
\mathbb{P} \left(\bigcap_{\eta} \bigcup_{k\geq 1} \bigcap_{N\geq k} \left\{\pi_{B}^{*}(N;\eta) \geq \frac{\lambda(N,\eta)}{2} \right\} \right)=1,
\end{align}
where $\eta$ is made to run over the fractions $(1/n)_{n \in \N}$. 

\begin{proof}
For fixed $N$ and $\eta$, write $M:= \pi_{\N}^{*}(N;\eta)$ and let $X \sim B(M , \lambda(N,\eta))$ denote a Binomial random variable ( i.e. a sum of $M$ independent Bernoulli trials, each of which yields success with probability $\rho^2$ ). Observe that for any natural number $k$ one has that
$$\mathbb{P}(\pi_{B}^{*}(N;\eta) \leq k)\leq \sum_{j=0}^{k} {{M}\choose{j}} (\rho^2)^{j} (1-\rho^2)^{M-j} = \mathbb{P}(X \leq k).$$
It now follows after a simple application of Chernoff's inequality (see for instance \cite[Theorem 1.8]{TaoVu} ) that
$$\mathbb{P}(X \leq \lambda/2) \leq \mathbb{P}(|X-\mathbb{E}(X)|\geq  \lambda/2) \leq e^{-\lambda}$$ 
and to conclude the proof of (\ref{heuristics}), we need only invoke the Borel-Cantelli lemma.
\end{proof}
\end{lemma}

Despite these heuristics one can piece together large subsets of the primes which have only large prime gaps. From \cite[Theorem 3]{PTIV} we get the existence of a constant $C>0$ such that
$$|\left\{ p_n \leq N | \  p_{n+1}-p_n \leq h   \right\}| \leq C \min \left( \frac{h}{\log N}, 1 \right) \pi(N)$$
for any pair of positive integers $N$ and $h$.
In particular, any hope of obtaining a result of the form ``If $\mathcal{C} \cap \mathcal{P}$ has positive relative density in the primes then $\mathcal{C}$ has small prime gaps" is dashed. Indeed, taking $h=\log N /2C$ in the statement just above, we gather that the collection of primes $p_n$ for which $p_{n+1}-p_{n}\geq 1/2C \log p_n$ constitutes a positive proportion of all prime numbers.
It thus becomes apparent that some kind of structure must be imposed on $\mathcal{C}$ if we wish to get primes in short intervals. We will explore two scenarios in which we are able to control the interaction of $\mathcal{C}$ with arithmetic progressions. In each setting we provide a somewhat general result and then give some examples of sets obeying the desired properties. The two main examples are

\begin{theorem}[Bohr sets have small prime gaps]\label{bohrsetshavesmallprimegaps}
Let $\mathfrak{g}(x)=\sum_{j=0}^D \alpha_j x^j \in \R[x]$ and suppose all coefficients $\alpha_j$ are Diophantine. Let $d \in (0,1)$, $\eta>0$ and write $\left\{x \right\} $ for the fractional part of a real number $x$. Then the Bohr set
$$\mathcal{A}:=\left\{ n \in \N \ | \  \left\{ \mathfrak{g}(n) \right\}  \in [0,d] \ \right\}$$
frequently contains $\eta$-small prime gaps and $H_{\mathcal{A}}(m) \ll_d e^{4m + \varepsilon}$ for $\varepsilon>0$ arbitrary.
\end{theorem}

\begin{theorem}[Shifted sets of square-free integers have small prime gaps]\label{kfreehavesmallprimegaps}
Let $\eta>0$ and $a \in \N$ be arbitrary. Then the set of shifted square-free integers 
\[ \mathcal{B}:=\left\{ n \in \N \ | \  n+a \text{ is square-free} \right\}\]
frequently contains $\eta$-small prime gaps. 
\end{theorem}

\begin{remark}
In Theorem \ref{bohrsetshavesmallprimegaps} one expects that the statement remains true as long as the coefficient $\alpha_D$ is irrational, however with our current methods this seems out of reach. We will discuss this in greater detail in section 6. In Theorem \ref{kfreehavesmallprimegaps} we do not obtain a bound for $H_{\mathcal{B}}(m)$. This is because our method relies on the ability to establish the correct order of magnitude $\pi_{\mathcal{B}}(N;\eta) \gg_{\eta} N/(\log N)$ and we are currently unable to do so for bounded prime gaps. 
\end{remark}
The case where $\mathcal{C}$ is an arithmetic progression has been studied by several authors. These results, which will be stated in section 6, differ from our own in the sense that we obtain small prime gaps frequently, as opposed to ``infinitely often". For other results concerning small gaps in special sets of primes we refer the reader to \cite{Cas},\cite{Pol}, \cite{Thorner}.
\\
\\
{\bf Notation}
We introduce some standard notation that will be used throughout the paper. For functions $f$ and $g$ we will use the symbols $f\ll g$ and $f=O(g)$ interchangeably to express Landau's big O symbol. A subscript of the form $\ll_{\eta}$ means the implied constant may depend on the quantity $\eta$. The statement $f \sim g$ means $f$ and $g$ are asymptotically equivalent, i.e., $\lim_{x \rightarrow \infty} f(x)/g(x)=1$. For two positive integers $l,m$ we write $l\asymp m$ when $l/2\leq m \leq 2l$ and the superscript $\flat$ indicates a summation over squarefree variables. We reserve the letter $\mu$ for the M\"obius function and $\Lambda$ for the von Mangoldt function. In place of $\Lambda$ we will sometimes use the function 
\begin{equation*}
\theta(n)= \left\{
\begin{array}{ll}
\log n & \mathrm{if}\ n \text{ is prime} \\
0 & \mathrm{otherwise.}
\end{array}
\right.
\end{equation*}
\section{The main results}
\subsection{Type A sets}
We first consider sets $\mathcal{A} \subset \Z$ which exhibit an even distribution among arithmetic progressions of any given modulus. 
Fix a natural number $k$. We say that a $k$-tuple $\mathcal{H}=(h_1,...,h_k)$ has height $h$ if each member is bounded in size by $h$. A set $\mathcal{A}$ will be of type $\operatorname{A}_k$ if it exhibits the following properties.

\begin{itemize}
\item[(a)] (Estimates for progressions in $\mathcal{A}$)
\\
There exist constants $c_1(k), c'>0, \theta>0 $  such that, for all positive integers $N$,

\[ \sum_{ \substack{n\leq N, n+\mathcal{H} \subset \mathcal{A} \\ n \equiv a \bmod q}} 1 = c_1(k) \frac{N}{q} +\mathcal{R}_{\mathcal{A}}(N,a,q) \]
and for any $\varepsilon>0$, the remainder term satisfies the bound
\begin{align}\label{averror1}
\sum_{q \leq Q} \max_{a } \left|\mathcal{R}_{\mathcal{A}}(N,a,q) \right| \ll_{\varepsilon} N^{1-c'} 
\end{align}
in the range $Q \leq N^{\theta}$ and uniformly over all $k$-tuples $\mathcal{H}$ of height $h \leq \log N$.

\item[(b)] (Estimates for primes in $\mathcal{A}$)
\\
There exists a constant $c_2(k)>0$ such that
\begin{align}\label{primesinA}
\sum_{ n\leq N, n+\mathcal{H} \subset \mathcal{A} }  \Lambda(n) = c_2(k) N +O\left( \frac{N}{(\log N)^C} \right)
\end{align}
for any $C>0$. In addition, for every pair $A,B>0 $, we have the estimate
\begin{align}\label{averror2}
\sum_{q \leq Q'} \sum_{\substack{\chi \bmod q\\ \chi \neq \chi_0}} \left| \sum_{ \substack{n\leq N\\ n+\mathcal{H} \subset \mathcal{A}} } \Lambda(n) \chi(n) \right| \ll \frac{N}{(\log N)^B}
\end{align}
uniformly in the range $Q'\leq (\log N)^A$. The index $\chi \bmod q$ runs over Dirichlet characters of modulus $q$ and $\chi_0$ denotes the trivial character.

\item[(c)](A bilinear form estimate for $\mathcal{A}$)
\\
\\
There exists a constant $C=C(\mathcal{A})$ such that
\begin{align}\label{bilinear1} 
\sum_{q \leq Q} \frac{q}{\varphi(q) }\sum_{\chi \bmod q}^{*} \left| \mathop{\sum_{m \leq M} \sum_{l\leq L}}_{ml + \mathcal{H} \subset \mathcal{A}} a_m b_l \chi(ml) \right| \ll  (M+Q^2)^{1/2}(L+Q^2)^{1/2}(\log ML)^C \left\|a\right\|_2 \left\|b\right\|_2 
\end{align}
for any pair of integers $M,L$ and $Q \leq \max(M^{\theta},L^{\theta})$. The sum $\sum_{\chi}^{*}$ indicates we are summing over primitive characters and $\left\|.\right\|_2$ denotes the $\ell^2$ norm of a sequence.
\end{itemize}
 
In this setting we will prove the following result.
\begin{theorem}\label{main1}
Any set $\mathcal{A} \subset \N$ of type $\operatorname{A}_k$ frequently contains $\eta$-small prime gaps, for $\eta>0$ arbitrary. Assuming that $c_3:=\sup_{l \in \N} c_1(l)/c_2(l) <\infty$, one has that 
$$H_{\mathcal{A}}(m) \ll \exp \left(2m/\theta +2c_3/ \theta +\varepsilon \right)$$
for all natural numbers $m$ and $\varepsilon>0$ arbitrary. 
\end{theorem}
With some additional structure on $\mathcal{A}$, we can obtain an improved result for type $A_1$ sets.

\begin{definition}
Let $\mathcal{G}$ be a finite set of primes and write $\mathcal{P}(\mathcal{G})=\prod_{p \in \mathcal{G}} p$. We say that a set $\mathcal{Y} \subset \N$ is locally well distributed if it satisfies estimates of the type
$$\sum_{ \substack{n\leq N \\ n \equiv a \bmod q}} 1_{\mathcal{Y}}(n) = \lambda \frac{N}{q} +\mathcal{R}_{\mathcal{Y}}(N,a,q)$$ 
for some $\lambda>0$ and there exists a constant $\epsilon>0$ such that 
$$\sum_{\substack{q \leq Q\\ (q,\mathcal{P}(\mathcal{G}))=1}} \max_{a} \left|\mathcal{R}_{\mathcal{Y}}(N,a,q) \right| \ll \frac{N}{(\log N)^A}$$
for $Q\leq N^{\epsilon}$ and $A>0$ arbitrary.
\end{definition}

\begin{definition}
Let $\mathcal{A}$ be a subset of the natural numbers for which we can write $\mathcal{A}' + \mathcal{Y} \subset \mathcal{A}$, where $\mathcal{A}'$ is a type $\operatorname{A}_1$ set and $\mathcal{Y}$ is locally well distributed. Then we say $\mathcal{A}=(\mathcal{A};\mathcal{A}', \mathcal{Y})$ is of type A'.
\end{definition}

\begin{Thm1prime}
Let $\mathcal{A}$ be a type $\operatorname{A}'$ set. Then for any $\eta>0$, $\mathcal{A}$ frequently contains $\eta$-small prime gaps. In addition one has the estimate $H_{\mathcal{A}}(m) \ll \exp \left(2m/\theta +2c_1(1)/(c_2(1) \theta) +\varepsilon \right)$ for all natural numbers $m$ $\varepsilon>0$ arbitrary. The constants $c_1(1),c_2(1)$ are associated to $\mathcal{A}'$. 
\end{Thm1prime}
We will apply these theorems to the case where $\mathcal{C}$ is a Bohr set. In section 6 we will prove that, after a somewhat careful selection of $k$-tuples $\mathcal{H}$, conditions $(a)-(c)$ hold. From this, Theorem \ref{bohrsetshavesmallprimegaps} will follow. For linear Bohr sets (i.e. when $D=1$) we will make use of Theorem \ref{main1}'. It should also be noted that, trivially, $\N$ is of type $A_k$ for any $k \in \N$, so that we recover the work of Goldston, Pintz and Y{\i}ld{\i}r{\i}m (e.g. \cite[Theorem 1]{PTIV}). 
\subsection{Type B sets}
We will call $\mathcal{B} \subset \Z$ a type B set if it has the form
$$\mathcal{B} =\bigcap_{m \in \mathfrak{M} } \bigcup_{a \in \mathfrak{M}_m} \left\{n \in \N| \ n \equiv a \bmod m \right\},$$
where $\mathfrak{M}$ is a collection of pairwise co-prime integers and associated to each $m \in \mathfrak{M}$ is a selection of residue classes $\mathfrak{M}_m \subset \Z_m $. We will write $\mathcal{N}_m$ for the complement of $\mathfrak{M}_m$ and assume that

$$\sum_{\substack{m \in \mathfrak{M}\\ m \geq x}}\frac{\mathcal{N}_m}{\varphi(m)} \ll \frac{1}{x^{\kappa}}.$$
for some $\kappa>0$ and any $x\geq 1$.
In this particular setting we will show that

\begin{theorem}\label{main2}
Let $\mathcal{B} \subset \N$ be a type B set satisfying the conditions listed above and let $\eta>0$ be arbitrary. Then $\mathcal{B}$ frequently contains $\eta$-small prime gaps. 
\end{theorem}

\begin{remark} Consider a type B set to which there are coupled infinitely many $m \in \mathfrak{M}$. To deal with sums of the form $\sum_{p} 1_{\mathcal{B}}(p)$, we are essentially counting primes in arithmetic progressions. However, since $\mathfrak{M}$ contains infinitely many moduli and we are not able to process a large quantity of remainder terms, venturing down this avenue would pose a problem. Instead we shall work with ``approximate type B sets $\mathcal{B}(z)$ ", which are collections of the form
\begin{align}\label{approxB}
\mathcal{B}(z):=\bigcap_{\substack{m \in \mathfrak{M}\\ m< z} } \bigcup_{a \in \mathfrak{M}_m} \left\{n \in \N| \ n \equiv a \bmod m \right\}.
\end{align}  
\end{remark}

At this point we note that Theorem \ref{kfreehavesmallprimegaps} follows easily from Theorem \ref{main2} but we will present the examples in section 6.
\\
\\
{\bf An outline of the strategy} 
\\
\\
To prove the main Theorems we will largely follow the framework set up in \cite{3}, \cite{PTIV} and \cite{May}. Let $\mathcal{H} \subset \N$ with $|\mathcal{H}|=k$ and for any such set, define the polynomial 
$$P_{\mathcal{H}}(n):=\prod_{h \in \mathcal{H}} (n+h).$$
Also, let $l< k$ and write $\mathcal{P}(z):=\prod_{p \leq z} p$. From now on $\mathcal{C}$ will denote either a type $\operatorname{A}_k$ set $\mathcal{A}$ or an approximate type B set $\mathcal{B}(z)$. In the latter case write $\mathcal{F}$ for the collection of primes dividing $\mathfrak{M}(z):=\prod_{m \in \mathfrak{M}, m \leq z} m$ and for type $\operatorname{A}_k$ sets we take $\mathcal{F}$ to be empty. When possible, we will treat both types of sets in a unified manner.\\
As in \cite{PTIV}, the first part of Theorem \ref{main1} is proven by evaluating weighted sums of the form

\begin{eqnarray}\label{THESIEVE}
\sum_{\mathcal{H}} \sum_{\substack {n=1 \\  (P_{\mathcal{H}} (n), \mathcal{P}(R^{\delta}) )=1}}^{N} \left( \sum_{\substack{h_0 \leq h\\n+h_0 \in \mathcal{C}}} \theta(n+h_0)-\log(3N) \right)  w(n)^2,
\end{eqnarray}
where $N$ is a natural number, $h:= \eta \log N$ and $\mathcal{H}$ runs over all $k$-tuples of height $h$. To obtain the estimates for $H_{\mathcal{A}}(m)$ we may drop the condition $(P_{\mathcal{H}} (n), \mathcal{P}(R^{\delta}) )=1$ and it becomes unnecessary to some over all $k$-tuples. It should be noted that the expression (\ref{THESIEVE}) can only be positive if there exists an interval of length $h$ in which $\mathcal{C}$ has at least two primes. In addition, (\ref{THESIEVE}) can only be ``large" if $\mathcal{C}$ has ``many" small prime gaps. 
Next we observe that, when splitting the above expression into two parts, the condition $n+h_0 \in \mathcal{C}$ falls only on the first sum. Intuitively, this indicates that an extra factor (the density of $\mathcal{C}$) will appear when dealing with this first sum but not for the second sum. To overcome this issue we will consider carefully selected $k$-tuples $\mathcal{H}$.

The weight function $w(n)$ was introduced in \cite{May} as a generalization of the GPY sieve

$$\Lambda_{R}(n,l):=\frac{1}{(k+l)!} \sum_{\substack{d| P_{\mathcal{H}}(n)\\ d\leq R} } \mu(d) \left(\log \frac{R}{d} \right)^{k+l}.$$

Let $\mathcal{H}:= \left\{ h_1, ..., h_k\right\}$ be an admissible $k$-tuple and suppose $F:\mathcal{R}_k \rightarrow \R$ is a differentiable function supported on the simplex $\mathcal{R}_k=\left\{ (x_1,...,x_k)\in [0,1]^k | \ \sum_{i=1}^k x_i \leq 1 \right\}$. Define 
$$\lambda_{d_1,d_2, ...,d_k}:=\left(\prod_{i=1}^k \mu(d_i) d_i \right) \sum_{\substack{r_1,...,r_k\\ d_i|r_i}}^{'} \frac{\mu(\prod_{i=1}^k r_i)^2 }{\prod_{i=1}^k \varphi(r_i)} F \left(\frac{\log r_1}{\log R}, ..., \frac{\log r_k}{\log R} \right),$$
where the superscript $'$ indicates the summation takes place over variables coprime to $W$. We will consider the sieve weights
$$w(n):=\sum_{\substack{d_i| n+h_i\\ (d_,\mathcal{F})=1 }} \lambda_{d_1,d_2, ...,d_k}$$ 
with corresponding sums 
\begin{align*}
S_1:= \sum_{\substack {n \sim N \\  (P_{\mathcal{H}} (n), \mathcal{P}(R^{\delta}) )=1}} w(n)^2, \ \ \ S_2:=\sum_{i=1}^k \sum_{\substack {n \sim N, n+h_i \in \mathcal{A} \\  (P_{\mathcal{H}} (n), \mathcal{P}(R^{\delta}) )=1}} \theta(n+h_i)  w(n)^2.
\end{align*}
To avoid the effect of small primes we let $W:=\prod_{p \leq D_0} p$ for some $D_0(k)$ and sum over values $n \equiv a_0 \bmod W$ for some $(a_0,W)=1$. We will also assume $\lambda_{d_1,...,d_k} $ is supported on variables coprime to $W$.
After establishing a Bombieri-Vinogradov type result in section \ref{sectionbombvin}, we will prove two crucial asymptotic estimates in section \ref{sectionasymp}, namely Propositions \ref{propasymp1} and \ref{propasymp2}, which are the analogues of \cite[Propositions 1 and 2]{PTIV}. In section \ref{sectioncompletingthms} these will be used to obtain an asymptotic formula for (\ref{THESIEVE}). On the other hand we will demonstrate a connection between the quantity (\ref{THESIEVE}) and  $\pi_\mathcal{C}(N;\eta)$, thereby proving our main results. 

\begin{definition}
Let $\mathcal{C}$ be a type $\operatorname{A}_k$ set $\mathcal{A}$ or an approximate type B set $\mathcal{B}(z)$. We introduce the constant
$$\mathfrak{S}_k(\mathcal{C}):=  \prod_{ \substack{p \in \mathcal{F}\\ D_0 < p <z }} \frac{1}{p}\left(1-\frac{1}{p} \right)^{k} \ \text{ when } \mathcal{C}=\mathcal{B}(z)$$
and $\mathfrak{S}_k(\mathcal{C})=1$ when $\mathcal{C}=\mathcal{A}$. For a piecewise differentiable function $F:\mathcal{R}_k \rightarrow \R$ we write
\begin{align*}
I_k:&= \int_{0}^{1}...\int_{0}^{1} F(t_1,...t_k)^2 dt_1,...dt_k,\\
J_k^{(m)}:&= \int_{0}^{1}...\left( \int_{0}^{1} F(t_1,...t_k) dt_m\right)^2 dt_1,...dt_{m-1} \ dt_{m+1}...dt_k
\end{align*}
\end{definition}

\begin{definition}[Selecting $k$-tuples]
Fix $k \in \N$ and let $\mathcal{H}=\left\{h_1,...,h_k \right\}$ be a $k$-tuple. We will write $\mathcal{H} \in \mathfrak{h}_k$ if, for each $1 \leq i \leq k$, we have $h_i \equiv 0 \bmod W$. 

\end{definition}

In the case of a  type $\operatorname{A}_k$ set we set $R:=N^{\theta/2 -\epsilon}$ for some $\epsilon >0$.

\begin{proposition}\label{propasymp1}
Fix $k \in \N$ and let $\mathcal{C}$ be a type $\operatorname{A}_k$ set or an approximate type B set $\mathcal{B}(z)$. Suppose $N^{c_1} \ll R^{\frac{1}{2+\delta}} \left(\log N \right)^{-c_2}$ for some constants $c_i$ depending on $k$. Furthermore, assume $\delta>0$ is sufficiently small compared to $k^{-1}$, $\mathcal{H} \in \mathfrak{h}_{k}$ is admissible and $h\ll \log R$ with $h \rightarrow \infty$ as $N \rightarrow \infty$. Then  we have
$$\sum_{\substack {n\leq N, n+\mathcal{H} \subset \mathcal{C} \\ (P_{\mathcal{H}} (n), \mathcal{P}(R^{\delta}) )=1}}^{\sharp} w(n)^2 = \gamma(\mathcal{C}, \mathcal{H}) \mathfrak{S}_k(\mathcal{C}) \frac{N (\log R)^k \varphi(W)^{k} }{W^{k+1}} I_{k}(F) (1+O(\delta k)).$$

When $\mathcal{C}=\mathcal{B}(z)$, the superscript $\sharp$ indicates we are selecting those values of $n$ satisfying $n+\mathcal{H} \subset \Z_m^{\times} \cap \mathfrak{M}_m$ for all $m \in \mathfrak{M}, m \leq z$ and $n \equiv a_0 \bmod W$. The constant factor is given by
\begin{equation*}
\gamma(\mathcal{C},\mathcal{H})=
\prod_{\substack{m \in \mathfrak{M}\\  D_0 < m \leq z}} \frac{\tau_m(\mathcal{H})}{m}  \ \ \ \text{  if  } \ \  \mathcal{C}=\mathcal{B}(z)
\end{equation*}
and $\gamma(\mathcal{C},\mathcal{H})$ is replaced by $c_1$ when  $\mathcal{C}=\mathcal{A}$. The symbol $\tau_m(\mathcal{H})$ represents the cardinality of the set $\cap_{i=1}^k \left( \left( \Z_m^{\times} \cap \mathcal{M}_m \right)-h_i \right)  $.
\end{proposition}

\begin{proposition}\label{propasymp2}
Let $\varepsilon>0$ be fixed. Given the same conditions as in Proposition 2, we have for $ N^{c_1} \ll R^{\frac{1/2 - \varepsilon}{2+\delta}} \left(\log N \right)^{-c_2}$ and $\mathcal{H}^0:=\mathcal{H} \cup \left\{h_0 \right\} \in \mathfrak{h}_{k+1}$ admissible

$$\sum_{\substack {n\leq N, n+\mathcal{H} \subset \mathcal{C} \\ (P_{\mathcal{H}} (n), \mathcal{P}(R^{\delta}) )=1}}^{\sharp} \theta(n+h_m) w(n)^2 = \gamma(\mathcal{C},\mathcal{H} ) \mathfrak{S}_k(\mathcal{C}) \frac{N (\log R)^{k+1} \varphi(W)^{k} }{W^{k+1}}  J_{k}^{(m)}(F) \left(1+O(\delta k+ G(\mathcal{H}^0) ) \right),$$

when  $h_m \in \mathcal{H}$ and in the case $h_0 \notin \mathcal{H}$, one has
$$\sum_{\substack {n\leq N, n+\mathcal{H} \subset \mathcal{C} \\ (P_{\mathcal{H}} (n), \mathcal{P}(R^{\delta}) )=1}}^{\sharp} \theta(n+h_0) w(n)^2 = \gamma(\mathcal{C},\mathcal{H}^0 ) \mathfrak{S}_{k-1}(\mathcal{C}) \frac{N (\log R)^{k} \varphi(W)^{k-1} }{W^{k}}  I_{k}(F) \left(1+O(\delta k+ G(\mathcal{H}^0) ) \right).$$

The remainder $G(\mathcal{H}^0)$ satisfies
$$\sum_{\substack{\mathcal{H}^0 \in \mathfrak{h}_{k+1}\\ h_i \ \text{distinct} }} G(\mathcal{H}^0) =o(h^{k+1}) \ \text{as} \  D_0 \rightarrow \infty.$$
\end{proposition}

\section{A Bombieri-Vinogradov theorem for type $\operatorname{A}_k$ sets}\label{sectionbombvin}

In this section we make the necessary changes to the proof of the Bombieri-Vinogradov theorem to ensure an identity of the form
$$\sum_{\substack{1 \leq m\leq N,  m+\mathcal{H} \subset \mathcal{A} \\  m \equiv a \bmod q } }  \theta(m)= \frac{c_2 N}{\varphi(q)} +E(N;q,a),$$
with the error term obeying the bound
\begin{align}\label{bombvin}
\sum_{q \leq N^{\theta-\varepsilon}} \max_{ \substack{ a \\ (a, q)=1}} |E(N;q,a)| \ll_{A,\varepsilon} \frac{N}{(\log N)^A}
\end{align}
for any $\varepsilon>0$ and $A>0$. In the next section we will need this result to demonstrate Propositions \ref{propasymp1} and \ref{propasymp2} for type A sets. 
\begin{remark}
To avoid any additional assumptions on type $A_k$ sets we wish to forgo the use of the P\'olya -Vinogradov theorem. This result is used in the proof of the Bombieri-Vinogradov theorem and if one were to follow the proof word-for-word in our current setting, one would require a bound of the form
$$\left| \sum_{n \leq N, n+ \mathcal{H} \subset \mathcal{A}} \chi(n) \right| \ll_{\varepsilon} N^{1/2}, $$
for any non-trivial character $\chi \bmod q$ with $1< q \leq N^{1/2-\varepsilon}$. For this reason, we will rearrange the argument in \cite[Chapter 28 ]{Dav} to better suit our purposes. 
\end{remark}

Let us begin by defining the sums
\begin{align*}
\psi_{\mathcal{A}}(N,\chi):=\sum_{\substack{1 \leq m\leq N\\ m+\mathcal{H} \subset \mathcal{A}  }} \Lambda(m) \chi(m) \ \text{ and} \ \ \  \psi_{\mathcal{A}}(N;q,a):=\sum_{\substack{1 \leq m\leq N \\ \substack{ m+\mathcal{H} \subset \mathcal{A} \\ m \equiv a \bmod q } } } \Lambda(m)=
\frac{1}{\varphi(q)}\sum_{\chi \bmod q} \sum_{\substack{1 \leq m\leq N\\  m+\mathcal{H} \subset \mathcal{A} } } \overline{\chi}(a) \Lambda(m) \chi(m). 
\end{align*}
Removing the expected main term from $\psi_{\mathcal{A}}(N;q,a)$, we may write

$$\psi_{\mathcal{A}}(N;q,a)-\frac{c^k N}{\varphi(q)}= \frac{1}{\varphi(q)} \sum_{\chi \bmod q}  \overline{\chi}(a) \psi_{\mathcal{A}}' (N,\chi) $$
where

\begin{equation*}
\psi_{\mathcal{A}}' (N,\chi)=\left\{
\begin{array}{ll}
\psi_{\mathcal{A}}(N, \chi) & \mathrm{when}\ \chi\neq \chi_{0}\\
\psi_{\mathcal{A}}(N, \chi_0) - c^kN & \mathrm{when}\ \chi =\chi_{0}.
\end{array}
\right.
\end{equation*}

For small moduli, e.g. $q \leq (\log N)^A$, we can estimate the sum

\[\sum_{q \leq Q} \frac{1}{\varphi(q)}  \sum_{\chi \bmod q}  \left| \psi_{\mathcal{A}}' (N,\chi) \right|\]
immediately, using the conditions imposed in ($b$). The resulting contribution is bounded by the RHS of (\ref{bombvin}).
\\
\\
{\bf An estimate for large moduli}
\\
\\For large values of $q \leq Q$ we will use the estimate

\begin{proposition}
For any $\varepsilon>0$ and $A>0$, there exists a constant $C>0$ such that 
\begin{align}\label{largemoduli}
\sum_{q \leq Q} \frac{q}{\varphi(q) } \max_{(a,q)=1}\left|\sum_{\chi \bmod q}   \overline{\chi}(a) \psi_{\mathcal{A}}'(N,\chi) \right| \ll \left( Q^{2} N^{1/2}+N^{5/6 }Q+N  \right)(\log QN)^C +\frac{QN}{(\log N)^A}
\end{align}
in the range $Q \leq N^{1/2- \varepsilon}$.
\end{proposition} 
{\bf Proof.} We will rework the argument laid out in \cite[chapter 28]{Dav} in conjuction with the inequality(\ref{bilinear1}). Once (\ref{largemoduli}) is established, the proof of  (\ref{bombvin}) is straightforward.\\

With this in mind, we decompose our sum, as described in \cite[p. 139]{Dav}, to obtain for $\chi \neq \chi_0$
$$\psi_{\mathcal{A}}'(N,\chi):=\sum_{\substack{m\leq N\\m+\mathcal{H} \subset \mathcal{A}} }   \Lambda(m)  \chi(m)=\sum_{m\leq N  }   \Lambda(m) f(m)=S_1+S_2+S_3+S_4$$
where 
$$S_1=\sum_{n \leq U} \Lambda(n)f(n)\ll U, \ \ \ \ \ \ S_2=-\sum_{t \leq UV} \left( \sum_{\substack{md=t\\m\leq U, d\leq V }} \mu(d)\Lambda(m)  \right) \sum_{r\leq N/t} f(rt),$$

$$S_3=\sum_{d \leq V} \mu(d) \sum_{h \leq N/d} f(dh) \log h \ \ \ \text{ and } \ \ \ \ \ S_4=-\sum_{U \leq m \leq N/V} \Lambda(m) \sum_{V<k\leq N/m } \left(\sum_{\substack{d|k\\d \leq V}} \mu(d) \right) f(mk).$$

From now on we will take $U=V=N^{\epsilon}$ for some constant $0<\epsilon<1$. To bound the sum $S_2$ we first restrict the variable $t$ to the range $[1,U]$ and then to the range $[U,UV]$. Denote the resulting sums as $S_2=S_2 '+ S_2 ''$. The sums $S_2 '$ and $S_3$ can be combined to create some degree of cancellation by noting that
\begin{align*}
\sum_{\chi \bmod q} \overline{\chi}(a) \left( S_2 ' + S_3\right)&=\sum_{t \leq U} \left( -(\mu \star \Lambda)(t) \sum_{r\leq N/t} \sum_{\chi \bmod q}\overline{\chi}(a) f(rt) + \mu(t) \sum_{r \leq N/t} \sum_{\chi \bmod q}\overline{\chi}(a) f(rt) \log r \right)\\
&= \varphi(q) \sum_{t \leq U} \left( \mu(t) \log t \sum_{ \substack{r\leq N/t \\ rt \equiv a \bmod q}}  f(rt) + \mu(t) \sum_{\substack{r\leq N/t \\ rt \equiv a \bmod q}} f(rt) \log r \right)\\
&= \varphi(q) \sum_{t \leq U} \mu(t) \sum_{\substack{r\leq N/t \\ rt \equiv a \bmod q}} f(rt) \log (rt)\\
&=\varphi(q) \sum_{\substack{t \leq U\\ (t,q)=1}}  \mu(t) \left( c^k\frac{N}{tq}\log N -c^k\frac{N}{tq}- \int_{1}^{N} \frac{\mathcal{R}(x,tq)}{x} \ dx \right).
\end{align*}
where $\star$ denotes a Dirichlet convolution and the last identity follows from summation by parts. For the sum $\sum_{\substack{\chi \bmod q\\ \chi \neq \chi_0}} \left( S_2 ' + S_3\right)$ we simply subtract the $\chi= \chi_0$ term in the first line. The additional term will be of the same form as the last line, just above, and can be dealt with in the same way. This last expression may be recast as
\begin{align}\label{mobius progression}
c\frac{\varphi(q)}{q} \sum _{l|q } \mu(l) \sum_{\substack{t \leq U\\ t \equiv 0 \bmod l}}  \frac{\mu(t)}{t} \left( N\log N -N \right)+ O\left( \int_{1}^{N} \varphi(q) \sum_{\substack{t \leq U\\ (t,q)=1}} \mathcal{R}(x,tq) \ \frac{dx}{x}  \right).
\end{align}
To estimate (\ref{mobius progression}) we require the following lemma. 
\begin{lemma}
Let $\Phi_{c}(x):=\exp((\log x)^c)$, then for any positive, squarefree integer $l \leq \Phi_{1/4}(x)$, one has the bound
$$\left|\sum_{\substack{t \leq x\\ t \equiv 0 \bmod l}} \frac{\mu(t)}{t} \right| \ll \frac{\omega(l)^{\omega(l)}}{\Phi_{1/4}(x)} ,$$
where $\omega(l)$ denotes the number of primes dividing $l$.
\begin{proof}
The bound is easily demonstrated by induction on $\omega(l)$. To be precise, we will show that 
\begin{align}\label{induction}
\left|\sum_{\substack{t \leq x\\ t \equiv 0 \bmod l}} \frac{\mu(t)}{t} \right| \ll (2\omega(l))^{\omega(l)} \exp \left(-\left( \log \left( \frac{x}{(\Phi_{1/4}(x))^{2\omega(l)}} \right) \right)^{1/4} \right).
\end{align}
When $\omega(l)=0$ and hence $l=1$, this estimate is obtained as a consequence of the bound $|\sum_{n \leq x} \mu(n)|\ll x/\Phi_{1/4}(x)$ (see for instance \cite[Chapter II.5]{Tenenbaum}) followed by partial summation. Let us now suppose that (\ref{induction}) has been established whenever $\omega(l')<r$ and let $l$ be such that $\omega(l)=r$. Furthermore, given $j\leq r$, let $l_j$ denote a generic divisor of $l$ having $j$ prime factors. Then we have
\begin{align*}
\sum_{\substack{t \leq x\\ t \equiv 0 \bmod l}} \frac{\mu(t)}{t}&= \frac{\mu(l)}{l} \sum_{\substack{t \leq x/l\\ (t,l)=1}} \frac{\mu(t)}{t}= \frac{\mu(l)}{l} \left( \sum_{j \leq r-1} (-1)^j {{r}\choose{j}} \sum_{\substack{t \leq x/l\\ t \equiv 0 \bmod l_j}} \frac{\mu(t)}{t} + \sum_{\substack{t \leq x/l\\ t \equiv 0 \bmod l}} \frac{\mu(t)}{t}\right)\\
&=\frac{\mu(l)}{l} \sum_{\substack{t \leq x/l\\ t \equiv 0 \bmod l}} \frac{\mu(t)}{t} +O\left( \frac{2^{r-1}}{l} \sum_{j \leq r-1} {{r}\choose{j}} (r-1)^{j} \exp \left[-\left( \log \left( \frac{x}{l(\Phi_{1/4}(x))^{2r-2}} \right) \right)^{1/4} \right] \right)\\
&=\frac{\mu(l)}{l} \sum_{\substack{t \leq x/l\\ t \equiv 0 \bmod l}} \frac{\mu(t)}{t} +O\left( \frac{2^{r-1}}{l}r^r \exp \left[-\left( \log \left( \frac{x}{l(\Phi_{1/4}(x))^{2r-2}} \right) \right)^{1/4} \right] \right).
\end{align*}
Repeating this process $K$ times yields
$$\left|\sum_{\substack{t \leq x\\ t \equiv 0 \bmod l}} \frac{\mu(t)}{t} \right| \ll \frac{1}{l^K} +2^{r-1} \omega(l)^{\omega(l)} \sum_{m \leq K} \frac{1}{l^m} \exp \left[-\left( \log \left( \frac{x}{l^m(\Phi_{1/4}(x))^{2r-2}} \right) \right)^{1/4} \right], $$
which gives the desired estimate after selecting the smallest value $K$ for which $l^K\geq \Phi_{1/4}(x)$. 
\end{proof}
\end{lemma}

In order to estimate the first sum in (\ref{mobius progression}) we apply the previous lemma for values of $l \leq \Phi_{1/8}(N)$ (combined with the fact that $\omega(l)\ll \log l$) and use a trivial bound otherwise. Together with (\ref{averror1}) it follows that
$$\sum_{q \leq Q} \frac{1}{\varphi(q)} \left| \sum_{\chi \bmod q} \left( S_2 ' + S_3\right) \right| \ll \frac{N}{(\log N)^{A}}$$

for any $A>0$. To estimate $S_4$ and $S_2''$ we will make use of the bilinear form inequality (\ref{bilinear1}). Since, currently, we are considering sums which run over complete sets of characters for each modulus $q$, these must first be converted to sums involving only primitive characters.\\
With this in mind, let $\chi \bmod q$ be generated by the primitive character $\chi_1 \bmod q' $ and write
$$a_m=\Lambda(m), \ \ \ \ b_k=\sum_{\substack{d|k\\d \leq V}} \mu(d).$$ 
First observe that
\begin{align*}
\sum_{\substack{U \leq m,k \leq N/U \\ mk \leq N, mk + \mathcal{H} \subset \mathcal{A}}} a_m b_k  \left[ (\chi_1-\chi)(mk) \right] =\sum_{\substack{U \leq m,k \leq N/U \\ mk \leq N}} a_m b_k  f_1(mk) - \sum_{\substack{U \leq m,k \leq N/U \\ mk \leq N, (mk,q)=1}} a_m b_k  f_1(mk).
\end{align*}

With regards to $S_4$, it follows that the contribution made by the character $\chi_1$ to the LHS of (\ref{largemoduli}) does not exceed
\begin{align*}
& \sum_{\substack{q \leq Q\\q \equiv 0 \bmod q' }  }^{\flat} \frac{1}{\varphi(q)}  \left|\sum_{\substack{U \leq m,k \leq N/U \\ mk \leq N, (mk,q)=1}} a_m b_k  f_1(mk) \right|=  \sum_{\substack{q \leq Q\\q \equiv 0 \bmod q'}  }^{\flat} \frac{1}{\varphi(q)}  \left| \sum_{l|q} \mu(l) \sum_{\substack{U \leq m,k \leq N/U \\ mk \leq N, l|mk}} a_m b_k  f_1(mk) \right| \\
& \ll \sum_{\substack{q \leq Q\\q \equiv 0 \bmod q' }  }^{\flat} \frac{1}{\varphi(q)} \sum_{l|\frac{q}{q'} } \left| \sum_{\substack{U \leq m,k \leq N/U \\ mk \leq N, l|mk}} a_m b_k  f_1(mk) \right| 
\ll \sum_{l \leq Q}^{\flat} \sum_{\substack{q \leq Q\\q \equiv 0 \bmod q' l}  }^{\flat} \frac{1}{\varphi(q)}  \left| \sum_{\substack{U \leq m,k \leq N/U \\ mk \leq N, l|mk}} a_m b_k  f_1(mk) \right|
\end{align*}
and hence we get that
$$\sum_{\substack{q \leq Q\\q \equiv 0 \bmod q' }  }^{\flat} \frac{1}{\varphi(q)}  \left|\sum_{\substack{U \leq m,k \leq N/U \\ mk \leq N, (mk,q)=1}} a_m b_k  f_1(mk) \right| \ll \sum_{l \leq Q}^{\flat}  \frac{\log N}{\varphi(l)} \frac{1}{\varphi(q')}  \left| \sum_{\substack{U \leq m,k \leq N/U \\ mk \leq N, l|mk}} a_m b_k  f_1(mk) \right|.$$
Summing over all primitive characters we conclude that $S_4$ accounts for a total bound of

\begin{align*}
 T_4:=&\sum_{l \leq Q}^{\flat}  \frac{\log N}{\varphi(l)} \sum_{q' \leq Q}^{\flat} \frac{1}{\varphi(q')}  \sum_{\chi}^{*} \left| \sum_{\substack{U \leq m,k \leq N/U \\ mk \leq N, l|mk}} a_m b_k  f_1(mk) \right|\\
  =& \sum_{l \leq Q}^{\flat}  \frac{\log N}{\varphi(l)} \sum_{q' \leq Q}^{\flat} \frac{1}{\varphi(q')} \sum_{\chi}^{*} \left| \sum_{\substack{U \leq m,k \leq N/U \\ mk \leq N}} a_m b_k  f_1(mk) \sum_{(d_1,d_2) \in P(l)} 1_{d_1,l}(m) 1_{d_2,l}(k)\right| \\
\ll & \sum_{l \leq Q}^{\flat}  \frac{\log N}{\varphi(l)} \sum_{(d_1,d_2) \in P(l)} \sum_{q' \leq Q}^{\flat} \frac{1}{\varphi(q')} \sum_{\chi}^{*} \left| \sum_{\substack{U \leq m,k \leq N/U \\ mk \leq N}} a_m 1_{d_1,l}(m) \ b_k 1_{d_2,l}(k) f_1(mk)   \right|.
\end{align*}
In the last two lines we have used $1_{d,l}(m)$ to indicate those integers $m$ for which $(l,m)=d$ and $P(l):=\left\{ (d_1,d_2) \subset \N^2 \left| \ d_1|l, \ d_2|l \text{ and } \ l|d_1d_2 \right. \right\}$. The three inner-most sums can now be dealt with as in \cite{Dav}, except that the estimate (\ref{bilinear1}) takes on the role of the large sieve. Indeed, by decomposing the range of $k$ into dyadic intervals $[M,2M]$ and setting $c_m:=a_m 1_{d_1,l}(m)$, $h_k:=  b_k 1_{d_2,l}(k)$ we find that on each such interval
\begin{align*}
\sum_{q \leq Q} \frac{q}{\varphi(q)} \sum_{\chi}^{*}  \left| \sum_{\substack{M \leq m \leq 2M \\ \substack{ m \in [U,N/U] \\ U \leq k \leq N/m }} }  c_m  h_k f_1(mk) \right| &\ll (M+Q^2)^{1/2}(\frac{N}{M}+Q^2)^{1/2}(\sum_{\substack{m \leq M}} \Lambda(m)^2)^{1/2} (\sum_{  k \leq N/M } d(k)^2)^{1/2} \log(N)^{C_1}\\
&\ll (M+Q^2)^{1/2}(\frac{N}{M}+Q^2)^{1/2}(M \log M)^{1/2}(\frac{N}{ M} (\log N)^3)^{1/2} \log(N)^{C_1}\\
& \ll (Q^2 N^{1/2}+\frac{QN}{M^{1/2}} + QN^{1/2}M^{1/2} + N) (\log N)^{C_1 +2}.
\end{align*}
  
The sum $S_2 ''$ may be treated as $S_4$ to obtain
\begin{align}
\sum_{q \leq Q} \frac{q}{\varphi(q)} \sum_{\chi}^{*} \left(|S_2 ''|+|S_4| \right) \ll  \left( Q^2 N^{1/2}+\frac{QN}{U^{1/2}} + QNU^{1/2} V^{1/2} + N \right)(\log N)^{C_2}.
\end{align}

Combining all of the above, we get 
$$\sum_{q \leq Q} \frac{q}{\varphi(q) }\sum_{\chi}^{*} \max_{y \leq N} \left| \psi_{\mathcal{A}}(x,\chi) \right| \ll \left( Q^2 N^{1/2}+N+\frac{QN}{U^{1/2}}+UQN^{1/2}+Q^{5/2}U  \right)(\log QUN)^{C_3} +\frac{QN}{(\log N)^A}$$
since $U=V$.
In the range $N^{1/3} \leq Q \leq N^{1/2}$ we choose $U=N^{2/3}/Q $. In this case it follows easily that all terms involving $U$ are bounded by
$$Q^2 N^{1/2}.$$

When $Q \leq N^{1/3}$ we select $U=N^{1/3}$. For such values of $Q$ we obtain the bound $QN^{5/6}$. This concludes the proof of (\ref{largemoduli}). To finish the proof of (\ref{bombvin}) we estimate the sum
$$\sum_{(\log N)^A \leq q \leq Q} \frac{1}{\varphi(q) } \max_{(a,q)=1}\left|\sum_{\chi \bmod q}   \overline{\chi}(a) \psi_{\mathcal{A}}'(N,\chi) \right|$$
by decomposing the range of $q$ into dyadic intervals and applying (\ref{largemoduli}) to each of the resulting sums.

\section{Obtaining the asymptotics}\label{sectionasymp}

Let $q=R^{\beta}$ denote a prime with $0 \leq \beta \leq \delta$. In order to prove Proposition \ref{propasymp1}, it is enough to demonstrate that 
$$ 
A_1:=\sum_{ \substack{ n\leq N\\ n+\mathcal{H} \subset \mathcal{C} }}^{\sharp}  w(n)^2 \sim \gamma(\mathcal{C})   \mathfrak{S}_{k}(\mathcal{C}) \frac{N (\log R)^k \varphi(W)^{k} }{W^{k+1}} I_{k}(F)
$$
and for each $1\leq l \leq k$,
$$A_1(q,l):=\sum_{\substack {n\sim N, n+\mathcal{H} \subset \mathcal{C}\\ \\ q|n+h_l } }^{\sharp} w(n)^2 \ll   \mathfrak{S}_{k}(\mathcal{C})\frac{N (\log R)^k \varphi(W)^{k} }{W^{k+1}} I_{k}(F) \left( \frac{\beta}{\varphi(q)}+ \frac{1}{\varphi(q)^2} \right).$$
To complete the argument one then sums the latter bound over all primes $D_0 <q \leq R^{\delta}$ .\\
These estimates are obtained by combining the ideas in \cite[Proposition 1]{PTIV} and \cite[Sections 4 and 5]{May} with some added details which we will point out.
Given an approximate type B set $\mathcal{B}(z)$, let $\mathfrak{M}(z)=\prod_{\substack{q \in \mathfrak{M}\\q <z }} q $ and take $\mathcal{F}$ to be the product of all primes dividing $\mathfrak{M}(z)$. The superscript $\sharp$, featuring in the definition of $A_1$, indicates we are selecting those values of $n$ satisfying $n+\mathcal{H} \subset \Z_m^{\times} \cap \mathfrak{M}_m$ for all $m \in \mathfrak{M}, m \leq z$.\\

{\bf Definition.} 
For $q=1$ or $q$ prime, $1 \leq l \leq k$ and $\mathcal{H}^{0}=\left\{h_0,h_1,...,h_k \right\}$, we introduce the quantities
\begin{align}
y_{r_1,...,r_k}(q,l)&=\left( \prod_{i=1}^k \mu(r_i) \varphi(r_i) \right) \sum_{\substack{d_1,...,d_k\\ \substack{r_i|d_i, (d_i,\mathcal{F})=1\\ q | d_l} }}  \frac{ \lambda_{d_1, ...,d_k} }{\prod_{i=1}^k d_i}, \\
x_{r_1,...,r_k}(q,l)&=\left( \prod_{i=1}^k \mu(r_i) g(r_i) \right) \sum_{\substack{d_1,...,d_k\\ \substack{r_i|d_i, (d_i,\mathcal{F})=1\\ q | d_l} }}  \frac{ \lambda_{d_1, ...,d_k} }{\prod_{i=1}^k \varphi(d_i)}
\end{align}  
and the analogous sums $y'_{r_1,...,r_k}(q,l)$, $x'_{r_1,...,r_k}(q,l)$ in which the condition $q| d_l$ is replaced by $q \nmid d_l$. Here, $g$ and $\rho$ are the totally multiplicative functions given by $g(p)=p-2$ and 
\begin{equation*}
\rho(p)= \left\{
\begin{array}{ll}
0 & \mathrm{if}\ p | \prod_{i=1}^k (h_0-h_i)  \\
1& \mathrm{ otherwise}\ .
\end{array}
\right.
\end{equation*}
on primes $p$. We also set $y_{\text{max}}:= \sup_{r_1,...,r_k} \left|y_{r_1,...,r_k} \right| $.

\begin{proposition}
As $N \rightarrow \infty$ one has the asymptotic formulas 
\begin{align}\label{A_1}
A_1(q,l)&=\frac{N}{W} \sum_{r_i} \frac{1 }{\prod_{i=1}^k \varphi(r_i)} \left( y_{r_1, ...,r_k}^2(q,l) + y_{r_1, ...,r_k}(q,l)y'_{r_1, ...,r_k}(q,l) + \frac{(y'_{r_1, ...,r_k})^2(q,l) }{q}\right)+ B_1(q,l)
\end{align}
\begin{align}\label{A_2}
A_2^{(i)}(q,l)&=\frac{N}{\varphi(W)} \sum_{r_i} \frac{1}{\prod_{i=1}^k g(r_i)} \left( x_{r_1, ...,r_k}^2(q,l) + x_{r_1, ...,r_k}(q,l)x'_{r_1, ...,r_k}(q,l) + \frac{(x'_{r_1, ...,r_k})^2(q,l) }{\varphi(q)}\right)  + B_2(q,l)
\end{align}
and the errors are bounded by
$$B_1(q,l)\ll  \frac{  y_{\text{max}}^2 N (\log R)^k}{\varphi(q)^2 W D_0}, \ \ B_2(q,l) \ll \frac{ N y_{\text{max}}^2}{\left( \log N \right)^A}+ \frac{  (  { y_{max}^{(i)}  } )^{2} N ( \varphi(W)^{k-2}) (\log R)^k}{ g(q)^2 W^{k-1} D_0}. $$
\begin{proof}
Expanding $A_1(q,l)$ we find that 
\begin{align}\label{A1asymp}
A_1(q,l)&=\sum_{\substack{d_1,...,d_k\\ e_1,..., e_k}} \lambda_{d_1, ...,d_k} \lambda_{e_1, ...,e_k} \sum_{\substack{n\sim N,n \in \mathcal{A}\\  \substack{q| n+h_l\\ [d_i,e_i] | n+h_i }}  }^{\sharp} 1 \\ \notag
&= \gamma(\mathcal{C})\frac{N}{W} \left( \sum_{\substack{d_1,...,d_k\\ \substack{e_1,..., e_k\\ q|[d_l,e_l]}}}^{*}  \frac{   \lambda_{d_1, ...,d_k} \lambda_{e_1, ...,e_k} }{ \prod_{i=1}^{k} [d_i,e_i]} + \sum_{\substack{d_1,...,d_k\\ \substack{e_1,..., e_k\\ q \nmid [d_l,e_l]}}}^{*}  \frac{   \lambda_{d_1, ...,d_k} \lambda_{e_1, ...,e_k} }{ q\prod_{i=1}^{k} [d_i,e_i]}\right) + \mathcal{R},
\end{align}
where the superscript $*$ indicates that $(d_i,\mathcal{F})=1$ and $(e_i,\mathcal{F})=1$ for all $i$ and we have used the important fact that $n+h_i \in \mathfrak{M} \cap \Z_m^{\times}$ for all $i$ to ensure that the constant $\gamma(\mathcal{C})$ appears. In the event that $q|M(z)$, clearly $A_1(q,l)$ vanishes and hence we may assume $q \nmid \mathcal{F}$. The error satisfies
$$\mathcal{R} \ll \sum_{\substack{d_1,...,d_k\\ e_1,..., e_k}} \left| \lambda_{d_1, ...,d_k} \lambda_{e_1, ...,e_k}\right| \ll \lambda_{max}^2 (\sum_{d<R} \tau_k(d))^2 \ll \lambda_{max}^2 R^2 (\log R)^{2k}.$$
The two sums on the RHS of (\ref{A1asymp}) may be separated into four parts, according to the divisibility of $d_l$ and $e_l$ by $q$. For example, following the manipulations in \cite[Lemma 5.1]{May} one arrives at   
\begin{align*}
\sum_{\substack{d_1,...,d_k\\ \substack{e_1,..., e_k\\ q\nmid d_l, q\nmid e_l}}}  \frac{   \lambda_{d_1, ...,d_k} \lambda_{e_1, ...,e_k} }{ \prod_{i=1}^{k} [d_i,e_i]}=&\sum_{u_1,...,u_k} \left( \prod_{i=1}^k \varphi(u_i)\right) \sum_{\substack{s_{i,j}\\ i \neq j }} \left( \prod_{\substack{1 \leq i,j \leq k\\ i \neq j}} \mu(s_{i,j})\right) \sum^{*}_{\substack{d_1,...,d_k\\ \substack{e_1,..., e_k\\\substack{ u_i| d_i, e_i\\ s_{i,j}| d_i, e_j \forall i\neq j}}}}  \frac{   \lambda_{d_1, ...,d_k} \lambda_{e_1, ...,e_k} }{ \prod_{i=1}^{k} (d_i e_i)}\\
&=\sum_{u_1,...,u_k} \left( \prod_{i=1}^k \frac{\mu(u_i)^2}{\varphi(u_i)}\right) \sum_{\substack{s_{i,j}\\ i \neq j }} \left( \prod_{\substack{1 \leq i,j \leq k\\ i \neq j}} \frac{\mu(s_{i,j}) }{\varphi(s_{i,j})^2}\right) y'_{a_1, ...,a_k}(q,l) y'_{b_1, ...,b_k}(q,l),
\end{align*}
where the superscript $*$ indicates the conditions $q\nmid d_l, q\nmid e_l$ and we have introduced the variables $a_j=u_j \prod_{i \neq j} s_{j,i}$ and $b_j=u_j \prod_{i \neq j} s_{i,j}$. Since $y'(q,l)$ is supported on variables which are coprime to $W$, we see that the main term in $A_1(q,l)$ comes from taking $s_{i,j}=1$ for all $i \neq j$. The remaining contributions must come from indices satisfying $s_{ij}>D_0$. These count towards an error no greater than
$$\frac{y^2_{max} N}{\varphi(q)^2 W} \left( \sum_{u\leq R, (u,W)=1}  \frac{\mu(u)^2 }{\varphi(u)}\right)^k \left( \sum_{s_{ij}>D_0 } \frac{\mu(s_{ij})^2 }{\varphi(s_{ij})^2}\right) \left( \sum_{t}  \frac{\mu(t)^2 }{\varphi(t)^{k^2-k-1}}\right) \ll \frac{y^2_{max} N \varphi(W)^{k} }{\varphi(q)^2 D_0 W^{k+1}}. $$
Applying the estimate $\lambda_{max} \ll y_{max} (\log R)^k$ (which was demonstrated in \cite[Lemma 5.1]{May}), the treatment of $A_1(q,l)$ is complete.\\
Turning to $A_2^{(0)}(q,l)$ we see that
\begin{align*}
&A_2^{(0)}(q,l)= \sum_{\substack{d_1,...,d_k\\ e_1,..., e_k}} \lambda_{d_1, ...,d_k} \lambda_{e_1, ...,e_k} \sum_{\substack{n\sim N, n \in \mathcal{A}\\  \substack{[d_i,e_i] | n+h_i\\ q|n+h_l }}  }^{\sharp} \theta(n+h_0)\\
&=\left( \prod_{\substack{m \in \mathfrak{M}\\ m > D_0}}  \frac{\tau_m(\mathcal{H}_0) }{ \varphi(m)} \right) \frac{N }{\varphi(W)} \left( \sum_{\substack{d_1,...,d_k\\ \substack{e_1,..., e_k\\q|[d_l, e_l] }}}^{*} \frac{  \rho(\prod_{i=1}^k [d_i,e_i]) \lambda_{d_i} \lambda_{e_i} }{ \prod_{i=1}^{k} \varphi([d_i,e_i])} + \sum_{\substack{d_1,...,d_k\\ \substack{e_1,..., e_k\\q \nmid [d_l, e_l] }}}^{*} \frac{  \rho(\prod_{i=1}^k [d_i,e_i]) \lambda_{d_i} \lambda_{e_i} }{ \varphi(q) \prod_{i=1}^{k} \varphi([d_i,e_i])} \right)+ \tilde{\mathcal{R} }.
\end{align*}

where we have used the shorthand $\lambda_{d_i}=\lambda_{d_1,..., d_k}$. The error term is no larger than

\begin{align*}
\tilde{\mathcal{R} } \ll \left(\prod_{m < z} \mathfrak{M}_m \right)  \ \lambda_{max}^2 \sum_{d \leq R^{2+\delta}}^{\flat}  \tau_k(d)  E'(d;N)& \ll_z y_{max}^{2} \left( \sum_{d \leq R^{2+\delta}}^{\flat} \frac{\tau_k(d)^2 N \log N}{\varphi(d)}  \right)^{1/2} \left(\sum_{d \leq R^{2+\delta}}^{\flat} E'(d;N) \right)^{1/2}\\
& \ll_z \frac{N y_{\text{max}}^2}{(\log N)^B}
\end{align*}
for any $B>0$. For the second inequality we used the trivial estimate $E'(N,t) \leq N (\log N) / \varphi(t) $. As with the discussion of $A_1(q,l)$ we split the two sums in the main term into four parts, depending on the divisibility of $d_l,e_l$ by $q$. The remainder of the proof goes through as in \cite[Lemma 5.2]{May}.
\end{proof}
\end{proposition}

To complete the proof of Propositions \ref{propasymp1} and \ref{propasymp2} we need the following sieve theory estimates.
\begin{lemma}\label{sieveestimates}
Let $\gamma$ be a multiplicative function satisfying 
$$0 \leq \frac{\gamma(p)}{p} \leq 1-A_1$$
$$-L\leq \sum_{w \leq p \leq z} \frac{\gamma(p)\log p}{p} - \kappa \log(z/w) \leq A_2$$
for all $2\leq w \leq z$. Let $g$ be the totally multiplicative function defined by $g(p)=\frac{\gamma(p)}{p- \gamma(p)}$ on primes and let $G:[0,1] \rightarrow \R$ be piecewise differentiable with $G_{max}:= \sup_{x \in [0,1]} (\left|G(x)\right| + \left|G'(x)\right|)$. Then for any prime $q=R^{\beta}$, we have the estimates
\begin{align}\label{sievesumq}
\sum_{\substack{d < z\\ q|d}}^{\flat}  g(d) G\left(\frac{\log d}{\log z} \right)= \frac{q-\gamma(q)}{q(q-1)} \mathfrak{S} \frac{(\log z)^{\kappa} }{\Gamma(k)} \int_{0}^{1} G(x)x^{\kappa-1} +O \left(\beta \frac{G_{max} }{\varphi(q)} (\log z)^{\kappa}  + \mathfrak{S} L G_{max} (\log z)^{\kappa-1} \right)
\end{align} 
and 
$$
\sum_{\substack{d < z\\ (q,d)=1}}^{\flat} g(d) G\left(\frac{\log d}{\log z} \right)= \frac{q-\gamma(q)}{q} \mathfrak{S} \frac{(\log z)^{\kappa}}{\Gamma(k)} \int_{0}^{1} G(x)x^{\kappa-1} +O \left(\beta \frac{G_{max} }{\varphi(q)} (\log z)^{\kappa} + \mathfrak{S} L G_{max} (\log z)^{\kappa-1} \right),$$ 
where the error terms depend only on $A_1,A_2,\kappa$ and
$$\mathfrak{S}=\prod_{p} \left(1-\frac{\gamma(p)}{p}\right)^{-1}  \left(1-\frac{1}{p} \right)^{\kappa}.$$
\begin{proof}
This is a rearrangement of \cite[Lemma 4]{GGPY}. The LHS of (\ref{sievesumq}) produces a main term with integral ranging in $[\log q/ \log z,1]$ rather than the full interval $[0,1]$. This accounts for the presence of the additional error $\frac{G_{max} \log q}{\log z}$ on the RHS of $(\ref{sievesumq})$.
\end{proof}

\end{lemma}

Choosing the smooth weights $y_{u_1,...,u_k}:=F \left(\frac{\log u_1}{\log R},...,\frac{\log u_k}{\log R} \right)$, we first observe that
\begin{equation*}
y_{u_1,...,u_k}(q,l)= \left\{
\begin{array}{ll}
y_{u_1,...,u_k} & \mathrm{if}\ q | u_l \\
-\frac{y_{u_1,...qu_l,...,u_k} }{\varphi(q)}& \mathrm{if}\ q \nmid u_l
\end{array}
\right.
\end{equation*}
and

\begin{equation*}
y'_{u_1,...,u_k}(q,l)= \left\{
\begin{array}{ll}
0 & \mathrm{if}\ q | u_l \\
y_{u_1,...,u_k}-\frac{y_{u_1,...qu_l,...,u_k} }{\varphi(q)} & \mathrm{if}\ q \nmid u_l.
\end{array}
\right.
\end{equation*}

It follows that 
$$A_1(q,l)=\frac{N}{W}\sum_{u_1,...,u_k} \frac{1 }{\prod_{i=1}^k \varphi(u_i)} \left( y_{u_1, ...,u_k}^2(q,l) + 2y_{u_1, ...,u_k}(q,l)y'_{u_1, ...,u_k}(q,l) +\frac{(y'_{u_1, ...,u_k})^2(q,l) }{q}\right)+E_1$$
and the sum in the main term becomes
$$\sum_{\substack{u_i\\ q|u_l}} \frac{y_{u_1, ...,u_k}^2 }{\prod_{i=1}^k \varphi(u_i)}  + \sum_{\substack{u_i\\ q \nmid u_l}} \frac{1 }{\prod_{i=1}^k \varphi(u_i)} \left( \frac{y_{u_1, ...,u_k}^2}{q} -\frac{2(1+1/q)}{\varphi(q)}y_{u_1, ...,u_k}y_{u_1,...qu_l,...,u_k} +y_{u_1,...qu_l,...,u_k}^2  \frac{3+1/q}{\varphi(q)^2}\right)$$
Applying Lemma \ref{sieveestimates} to the summation over each variable $u_i$ we conclude the discussion of $A_1(q,l)$.\\
Before carrying out the estimate for $A_2(q,l)$ we require a combinatorial rearrangement of the functions $x(q,l)$ and $x'(q,l)$.
\begin{lemma}
For a given $k+1$- tuple $\mathcal{H}^0= \left\{h_0, h_1,...,h_k\right\} $ with $h_0=h_m$ one has that
$$x_{r_1,...,r_k}(q,l)= \sum_{a_m} \frac{y_{r_1,...,a_m,...r_k}(q,l)}{\varphi(a_m)} + 
O\left(\frac{y_{max} \varphi(W)(\log R)}{W D_0}\right)$$
and the corresponding formula for $x'_{r_1,...,r_k}(q,l)$ is obtained by replacing $y_{r_1,...,a_m,...r_k}(q,l)$ with $y'_{r_1,...,a_m,...r_k}(q,l)$. Furthermore, when $h_0 \neq h_m$ for all $1 \leq m \leq k$, one has that
$$ x_{r_1,...,r_k}(q,l)= h^{k+1} y_{r_1,...,r_m,...r_k}(q,l) + O\left( \frac{y_{max} \varphi(W) (\log R)}{W D_0} B(\mathcal{H}^{0}) \right),$$
with the error $B(\mathcal{H}^{0})$ satisfying

\begin{equation*}
\sum_{\substack{h_0,h_1,...,h_k\\ \text{distinct} }} B(\mathcal{H}^{0}) \ll \left\{
\begin{array}{ll}
 h^{k+1}/\varphi(q)  & \mathrm{if}\ q | r_l \\
h^{k+1}  & \mathrm{if}\ q \nmid r_l.
\end{array}
\right.
\end{equation*}

\begin{proof}
We will focus on the second identity since the first is derived in much the same manner. For the same reason we will only discuss the case where $q|r_l$. Let us begin by inserting the identity
$$\lambda_{d_1,...,d_k}= \left( \prod_{i=1}^k \mu(d_i) d_i\right) \sum_{\substack{r_1,...,r_k\\ d_i|r_i }} \frac{y_{r_1,...,r_k}}{\prod_{i=1}^k \varphi(r_i)}$$
into the definition of $x_{r_1,...,r_k}(q,l)$ to get
$$x_{r_1,...,r_k}(q,l)=\left( \prod_{i=1}^k \mu(r_i) g(r_i)\right) \sum_{\substack{d_1,...,d_k\\ r_i| d_i}}\left( \prod_{i=1}^k \frac{\rho(d_i) \mu(d_i) d_i}{\varphi(d_i)}\right) \sum_{\substack{a_1,...,a_k\\ d_i|a_i }} \frac{y_{a_1,...,a_k}}{\prod_{i=1}^k \varphi(a_i)}.$$
Setting $s_i:=(a_i, \mathcal{F})$ and interchanging the order of summation we first note that
$$\sum_{\substack{d_i\\ r_i| d_i,  d_i| a_i}} \frac{\rho(d_i) \mu(d_i) d_i}{\varphi(d_i)}=\frac{\mu(r_i) \rho(r_i) r_i}{\varphi(r)} \prod_{ \substack{p| a_i \\  p\nmid \ r_i, s_i}} \left(1-\frac{p \ \rho(p)}{p-1}\right)=\frac{\mu(r_i) \rho(r_i) r_i h(a_i)}{\varphi(r_i) h(r_i) h(s_i)}.$$
Here, $h$ is the multiplicative function for which $h(p)=1-\rho(p) p/(p-1)$ on primes. We gather that
$$x_{r_1,...,r_k}(q,l)=\left( \prod_{i=1}^k \frac{r_i g(r_i) \rho(r_i)}{\varphi(r_i) h(r_i)} \right) \sum_{\substack{a_1,...,a_k\\ r_i|a_i }} y_{a_1,...,a_k} \prod_{i=1}^k\frac{h(a_i/s_i)}{ \varphi(a_i)}.$$
Taking into consideration the support of $y$ we see that the only non vanishing terms occur when $a_i=r_i$ or $a_i > r_i D_0$. Let us examine the latter situation. Suppose $a_j>r_j D_0$, and for each $1 \leq l \leq k$ write $t_l$ for the largest divisor of $a_l/s_l$ satisfying $\rho(t_l)=0$. Then we have 
$$\sum_{\substack{h_0,h_1,...,h_k\\ \text{distinct} }} \sum_{\substack{a_1,...,a_k\\ r_i|a_i }} \frac{h(a_i)}{ \varphi(a_i)} \ll_k \sum_{\substack{h_0,h_1,...,h_k\\ \text{distinct} }}\sum_{\substack{a_1,...,...,a_k\\ \substack{ t_1...t_k|h_1-h_0}\\a_j>r_j D_0 } } \frac{1}{ \varphi(a_l/t_l s_l)^2 \varphi(t_l)}.$$
Switching the order of summation in the last expression and recalling that $q$ divides $r_l$ (and hence $a_l$) we get the desired error term. The main term becomes $\prod_{i=1}^{k} (g(r_i) r_i \rho(r_i)/ \varphi(r_i)^2) y_{r_1,...,r_m,...r_k}(q,l) $. After noting that $g(p)p/\varphi(p)^2=1+O(p^{-2})$ and summing over all $k+1$- tuples $\mathcal{H}$ the result follows.
\end{proof}
\end{lemma}
The derivation of the estimate for $A_2(q,l)$ is similar to our discussion of $A_1(q,l)$. When $h_0=h_m$, equation (\ref{A_2}) combined with the previous lemma leads to the sum
$$\sum_{\substack{r_i\\i \neq m}} \frac{1}{\prod_{i=1}^k g(r_i)} \left( \sum_{a_m} \frac{y_{r_1,...,a_m,...r_k}(q,l)}{\varphi(a_m)} \right)^2.$$ 
We note that a factor $\left( \frac{\varphi(W)}{W}\prod_{\substack{p \in \mathcal{F}\\ p > D_0}}  \frac{\varphi(p) }{ p} \right)^2$ is introduced after applying Lemma \ref{sieveestimates} to the expression in brackets. Together with the factor appearing in equation (\ref{A_2}), this accounts for the constant $\gamma(\mathcal{C}, \mathcal{H})$ in Proposition \ref{main2}. The summations over the remaining variables are carried out as in \cite[Lemma 6.3]{May}.

\section{Completing the proofs of the main theorems}\label{sectioncompletingthms}

To complete the proofs of our main theorems we use the following key result.

\begin{proposition}
Define the quantity 
$$M_k:=\sup_{F} \frac{\sum_{m=1}^k J_k^{(m) }(F) }{I_k(F)},$$
where the supremum is taken over all differentiable functions $F:\mathcal{R}_k \rightarrow \R$ supported on the simplex $\mathcal{R}_k=\left\{ (x_1,...,x_k)\in [0,1]^k | \ \sum_{i=1}^k x_i \leq 1 \right\}$. Then $M_k > \log k -2 \log \log k -2$ for sufficiently large $k$.
\begin{proof}
This is \cite[Proposition 4.3]{May}. 
\end{proof}
\end{proposition}

To prove the first part of Theorems \ref{main1}, \ref{main1}' and Theorem \ref{main2}, we recall that \[ \pi_\mathcal{C}(2N;\eta)=|\left\{q_n \leq 2N | \  q_{n+1}-q_n \leq \eta \log N  \right\}|\]
and consider the sum

\begin{align*}
\tilde{S}_{R,\mathcal{C}}:= \frac{1}{h (\log R)^k} \sum_{\substack{N \leq n \leq 2N\\ n \in \mathcal{A}}} \left( \Theta_{\mathcal{C}}(n,h) - \log(3N) \right) \sum_{\mathcal{H}}^{\square} w(n)^2, 
\end{align*}
where
$$\Theta_{\mathcal{C}}(n,h)=\sum_{\substack{h_0 \leq h\\ n+h_0 \in \mathcal{C}}} \theta(n+h_0)$$
and the superscript $\square$ indicates that the summation takes place over admissible $k$-tuples $\mathcal{H} \subset [1,h]$ for which $(P_\mathcal{H}(n), \mathcal{P}(R^\delta) )=1$ and $\mathcal{H} \in \mathfrak{h}_{k+1}$. In the case of a type $\operatorname{A}'$ set, we will also assume that each member of $\mathcal{H}$ lives in $\mathcal{Y}$.
The upper bound for $\tilde{S}_{R,\mathcal{A}}$ is obtained in precisely the same manner as \cite{PTIV} with the addendum that
$$\sum_{\substack{N \leq n \leq 2N\\ \Theta(n,h)\geq 3/2 \log N}}  1_{\mathcal{A}'}(n) \leq h \ \pi_\mathcal{A}(2N;\eta) +O(N\exp(-c\sqrt{\log N})).$$
Observe that when $(P_\mathcal{H}(n), \mathcal{P}(R^\delta) )=1$,
$$w(n) \ll 2^{\frac{k\log R}{\delta \log N}} (\log R)^k. $$
It follows, as in \cite{PTIV}, that
\begin{align*}
 \tilde{S}_{R,\mathcal{C}} & \ll (h \ \pi_\mathcal{C}(2N;\eta) +o(\pi(2N)))^{1/2} \notag \\
 &\times \frac{(\log R)^{k} \log N}{N^{1/2}h^k} (2k+2)! 2^{(2k \log 3N)/\delta \log R} \left(\frac{h}{\delta \log R} \right)^{k/2} \left(1+\frac{h}{\delta \log R} \right)^{(k+2)/2}.
\end{align*}
On the other hand, combining the asymptotics with Proposition 5, one finds
\begin{align*}
\tilde{S}_{R,\mathcal{C}} \gg c_2 \mathfrak{S}_{k}(\mathcal{C})\frac{ {\varphi(W)^{k}} }{ W^{k+1}} \log N I_k(F)  \left( M_k \frac{\log R}{\log N} + c_2 \eta\frac{W}{\varphi(W) D_0}  -\frac{c_1}{c_2} + O(k^2 \delta)+O\left(\frac{1}{D_0(k)} \right) \right).
\end{align*}
In the case of a type $\operatorname{A}'$ set one gets a constant of the form $ c_{3}(\mathcal{G}, \lambda) c \eta$ instead of $c \eta $. By choosing $\delta \ll \eta \theta $ we get $1+h/(\delta \log R) \ll h/(\delta \log R)$. For $\eta$ sufficiently sufficiently small and $k$ sufficiently large, it is at once clear that $\tilde{S}_R >0$ and a small calculation shows that 

$$\pi_\mathcal{A}(2N;\eta) \gg C_{\mathcal{C} }(k,z)\eta^{r(k)}\pi(N),$$
for some positive integer $r(k)$. To finish the proof of Theorem \ref{main2} we now let $z:=C_{\mathcal{C} }(k)^{1/(\kappa +1)}$.Now set $\mathcal{D}:=\mathcal{B}(z)^{c} \cap \mathcal{B}$ and apply the Brun-Titchmarsh Theorem (see for example \cite[Theorem 6.6]{6}) to obtain
$$\pi_{\mathcal{D}}(N;\eta)\ll \pi_{\mathcal{D}}(N)\ll \left( \sum_{\substack{q \in \mathfrak{M}\\q \geq z}} \frac{\mathcal{N}_q}{\varphi(q)} \right) \pi(N)\ll \frac{\pi(N)}{z^{\kappa}}.$$\\
For the second part of Theorem \ref{main1} we note that
$$\sum_{\substack{N \leq n \leq 2N\\ n \in \mathcal{A}}} \left( \sum_{i=1}^k \theta(n+h_i) - \log(3N) \right) w(n)^2 \gg 
c_2 \mathfrak{S}_{k}(\mathcal{A}) \frac{ {\varphi(W)^{k}} }{ W^{k+1}} \log N I_k(F)  \left( M_k \frac{\log R}{\log N}  -\frac{c_1}{c_2} \right).$$
In order to make this sum positive we need only take $k \gg \exp \left(2m/\theta +2c_3/ \theta +\varepsilon \right)$.

\section{Some examples}
\subsection{Arithmetic progressions}
As a first example, we consider arithmetic progressions $\left\{ n \in \N \ | \ n \equiv a \bmod q \right\}$ where $(a,q)=1$. It is at once clear that the required conditions for a type B set hold and hence we get the following corollary.

\begin{proposition}
Let $\eta>0$ be arbitrary. Then any arithmetic progression $(nq+a)_{n \in \N}$ with $(a,q)=1$ frequently contains $\eta$-small prime gaps.

\end{proposition}

Several authors have studied the case of arithmetic progressions and obtained results similar to the above. Firstly, Goldston, Pintz and Y{\i}ld{\i}r{\i}m showed in \cite{PTIII} that there are small prime gaps in the progression $\left\{ n \leq N \ | \ n \equiv a \bmod q \right\}$ and one can even let $q$ grow slowly with $N$.

\begin{theorem}
Let
$\varepsilon$ and $A$ be arbitrary fixed positive numbers. Let $q$ and $N$ be arbitrary, sufficiently large integers, satisfying
$$q_0(A, \epsilon)\ll q \ll (\log \log N)^A, \ \ N> N_0(A, \epsilon), $$
and let
$a$ be arbitrary with $(a,q) = 1$. Then there exist primes $p,p' \in [N/3,N]$ such that $p' \equiv p \equiv a \bmod q$ and $p'-p < \varepsilon \log p$.
\end{theorem}

T. Freiberg demonstrated, in \cite{Frei}, that it is possible to find consecutive primes in short intervals which are both congruent to $a \bmod q$.  

\begin{theorem}
Fix any positive number $\epsilon$, and fix a pair of coprime integers $q \ge 3$ and $a$. There is an absolute positive constant $c$ such that, for all sufficiently large $X$,
\begin{align*}
\sum_{\substack{ p_{r} \le X, p_{r+1} - p_{r} < \epsilon \log p_{r}\\ p_{r} \equiv p_{r+1} \equiv a \bmod q} } 1 \geq X^{1 - c/\log\log X}.
\end{align*}
\end{theorem}

\subsection{Shifted sets of $k$-free numbers}
Our second application of Theorem \ref{main2} pertains to shifted sets of $k$-free numbers, i.e sets of the form
$$\mathcal{B}:=\bigcap_{p} \left\{n \in \N | \   n \not \equiv a \mod p^k \right\},$$
where $a$ is any fixed integer. Again, it is easily verified that $\mathcal{B}$ is a type B set so we obtain Theorem \ref{kfreehavesmallprimegaps}.

\subsection{Bohr sets}

Sequences of the type $(\left\{ \mathfrak{g}(n) \right\})_{n \in \N}$, with $\mathfrak{g}$ a polynomial, have been the subject of much study. It was demonstrated by H.Weyl that they are uniformly equidistributed in the unit interval, provided that the leading coefficient of $\mathfrak{g}$ is irrational. Later on, I.M. Vinogradov showed that the sequence remains equidistributed if one restricts $n$ to prime values. In this section we will add yet another result to this subject by proving Theorem \ref{bohrsetshavesmallprimegaps}. Let $\mathfrak{g}(x)=\sum_{j=1}^D \alpha_j x^j \in \R[x]$. A Bohr set is a collection of the form $\mathcal{A}:=\left\{ n \in \N \ | \  \left\{ \mathfrak{g}(n) \right\}  \in [0,d] \ \right\}$.

\begin{definition}
An irrational number $\alpha$ is of type $\rho>0$ if 
$$\rho = \sup \left\{ \gamma \in \R \ | \ \ \liminf_{m \rightarrow \infty} m^{\gamma} \left\| m\alpha \right\|=0  \ \right\}.$$
Here $\left\|x\right\|$ denotes the distance from $x$ to the nearest integer. A number $\alpha$ which obeys such a bound is said to be Diophantine.
\end{definition}

\begin{remark} Observe that Theorem \ref{bohrsetshavesmallprimegaps} holds for Lebesgue-almost all $D$-tuples $(\alpha_0,...\alpha_D) \in \R^D$ and as a consequence of the Thue-Siegel-Roth Theorem (see \cite{DavRoth}), whenever the $\alpha_j$ are algebraic irrationals.\\
We would like to show that $\mathcal{A}$ is a type $\operatorname{A}_k$ set (for large $k$) but there are some immediate algebraic obstructions that must be overcome. Consider, for example, the polynomial $\mathfrak{g}(x)= \sqrt{2}x^2 $. Observe that the events $n \in \mathcal{A}$, $n+1 \in \mathcal{A}$,$n+2 \in \mathcal{A}$ and $n+3 \in \mathcal{A}$ are not independent, since 
$$\left\{ \sqrt{2}n^2 \right\}-3\left\{ \sqrt{2}(n+1)^2 \right\}+3\left\{ \sqrt{2}(n+2)^2 \right\}-\left\{ \sqrt{2}(n+3)^2 \right\} \equiv 0 \bmod 1.$$
In other words, (\ref{averror1}) does not even hold for $4$-tuples.
\end{remark}

\subsubsection{ {\bf Some background information and tools }}
We begin by recalling some useful facts from the theory of Diophantine approximation (which can be found in \cite{Rav}). For the remainder of this section we will assume $\alpha$ is Diophantine of type $\rho$. 
Given any positive integer $M$, the collection $\left\{ m \leq M \ | \ \left\{m \alpha \right\} \in [0,d] \ \right\}$ has a very neat combinatorial structure. In order to describe this structure, rearrange the natural numbers up to $M$ in such a way that $\left\{s_1 \alpha \right\}<\left\{s_2 \alpha \right\}<...<\left\{s_M \alpha \right\}$. One has the recurrence relationship 
  
\begin{equation*}
s_{j+1}=\left\{
\begin{array}{lll}
s_j + s_1 & \mathrm{when}\ s_j \leq M-s_1\\
s_j+ s_1-s_M & \mathrm{when}\ M-s_1 < s_j < s_M\\
s_j-s_M	& \mathrm{when}\ s_M \leq s_j. 
\end{array}
\right.
\end{equation*}

Without loss of generality we may assume $||s_1 \alpha|| \leq ||s_M \alpha||$ so that $||s_M \alpha||, ||(s_1-s_M) \alpha||\geq 1/4M$. By Dirichlet's theorem we have that $||s_1 \alpha||\leq 1/M $ and since $\alpha$ is Diophantine, we easily find that $s_1 \gg_{\varepsilon} M^{1/\rho - \varepsilon}$ for any small $\varepsilon>0$. It also follows easily from the above that
\begin{align}\label{diosolns}
\left\{ b \leq M | \ \ \|b \alpha \| \leq \frac{1}{M^{1+c}}  \right\} \subset \left\{ b \leq M | \ b \equiv 0 \bmod s_1  \right\}
\end{align}
for any $c>0$ (and $M$ sufficiently large). Another important notion related to our problem is that of discrepancy. Given a sequence ${\bf x}=(x_1,...,x_N) $ in $[0,1]$ and real numbers $0\leq \delta <\beta \leq 1$,  write 
$$A([\delta,\beta)) = \left| \left\{ n \leq N \ | \ x_n \in [\delta,\beta) \right\} \right|.$$ 
The discrepancy of ${\bf x}$ is defined to be 
$$D_N(\omega)= \sup_{0\leq \delta <\beta \leq 1} \left|  \frac{A([\delta,\beta);N)}{N}-(\beta-\delta)  \right|.$$

\begin{remark}\label{fourierremark} Instead of working with the discontinuous characteristic function $\chi_{[0,d]} $ we shall sometimes make use of a smooth cut-off function $\psi(x): [0,1] \rightarrow [0,1]$ satisfying 

\begin{equation*}
\psi(x)=\left\{
\begin{array}{ll}
1 & \mathrm{when}\ \delta \leq x \leq d- \delta\\
0 & \mathrm{when}\ x \notin [0,d].
\end{array}
\right.
\end{equation*}
Employing such a bump function will improve the rate of convergence of $\psi$'s Fourier expansion. More precisely, if $\psi$ is $r$-times continuously differentiable with $\psi^{i}(0)=\psi^{i}(d)=0$ for $i<k$, one easily shows that the Fourier coefficients $c_k$ grow like

$$c_k \ll \min \left(\frac{1}{k}, \frac{N^{rC}}{k^{r}} \right)$$
provided that we choose $\delta:= N^{-C}$. In other words, for such a choice of $\delta$ we get rapid convergence as soon as $k \gg N^C$.
\end{remark}

The following results will aid us in the verification of conditions $(a)-(c)$. For more details, we refer the reader to \cite{KN}, \cite{Nath}, \cite{Har} and \cite{Stein} respectively. 
\begin{theorem}[Erd\"os-Tur\'an-Koksma]
Let ${\bf x}=(x_1,...,x_M)$ be a sequence of real numbers in the unit interval $[0,1]$. Then for any $H \in \N$ the discrepancy of this sequence is bounded by
$$D_M({\bf x}) \ll \frac{1}{H} + \sum_{r \geq 1} \frac{1}{r} \left| \frac{1}{n} \sum_{n=1}^{M} e(r x_n) \right|. $$
\end{theorem}

\begin{theorem}[Weyl's inequality]
Let $f(x) \in \R[x]$ be a polynomial of degree $k$ with leading coefficient $\alpha$ satisfying $|\alpha-r/\nu|\leq 1/\nu^2$ for some pair of coprime integers $r, \nu$ with $\nu>0 $. Then 
$$\left|\sum_{n \leq N} e(f(n)) \right| \ll N^{1+\varepsilon} \left(\frac{1}{\nu} + \frac{1}{N} +\frac{\nu}{N^k} \right)^{2^{1-k}}$$
for any $\varepsilon>0$.
\end{theorem}

%Before proceeding to the proof of Lemma \ref{smallmodBohr}, we recall the following estimate for exponential sums involving primes (see \cite{Har}).%
\begin{theorem}\label{Harman}
Let $f \in \R[x]$ be a polynomial of degree $k$, with leading coefficient $\alpha$ and suppose $|\alpha-a/v|\leq 1/v^2$ with $(a,v)=1$ and $v\leq N$. Then for any $\varepsilon >0$,
$$\left|\sum_{n=1}^{N}  \Lambda(n) e(f(n)) \right| \ll N^{1+\varepsilon}\left( \frac{1}{N^{1/2}} +\frac{1}{\nu} + \frac{\nu}{N^k} \right)^{4^{1-k}}.$$
\end{theorem}

In addition we will use a special case of van der Corput's lemma.
\begin{lemma}[van der Corput, special case]
Let $\phi(x) \in C^1([a,b])$ such that $|\phi'(x)|\geq 1$ and $\phi'(x)$ is monotonic in $(a,b)$. Then
$$\left|\int_{a}^{b} e(\phi(x)) \right| =O(1), $$
with the implied constant being independent of $a$ and $b$. 
\end{lemma}

\begin{lemma}
For $\alpha$ Diophantine and $\varepsilon>0$ fixed, one has the inequality 
$$ |\mathcal{B}_{\alpha}(M,M^{k+C},m)|=\left| \left\{ b \leq M | \ \ \|m b^k \alpha \| \leq 1/M^{k+C}  \right\} \right| \ll_{k, \varepsilon} M^{1-c(C, \rho)}$$ 
for any $C>0$  and uniformly in the range $1 \leq m \leq M^{C-\varepsilon}$. %Moreover, when $k=1$, one has $|\mathcal{B}_{\alpha}(M,M^{1-\varepsilon},m)|\ll M^{1-c(C, \rho)}$ uniformly in the range $1 \leq m \leq M^{\varepsilon}$. 
\begin{proof}
Let $\wp^k_{m} \subset \N$ denote the set of all integers of the form $m n^k$ and write $H:=\left\lceil N^{C}\right\rceil$. According to (\ref{diosolns}) there exists an integer $s_1 \leq M^k H$ such that
$$ \left|\mathcal{B}_{\alpha}(M,M^{k+C},h) \right| \leq \left|\left\{ ls_1| \ l \leq M^{k} H/s_1 \  \right\} \cap \wp^k_{h} \right|.$$
Since $\alpha$ is Diophantine, we have a lower bound of the type $s_1\geq (MH)^{c(\rho)}$. Combining this fact with the above estimate, the result easily follows.
\end{proof}
\end{lemma}

\subsubsection{ {\bf Verifying conditions $(a)-(c)$ for Bohr sets }}

To begin with, we will need to choose a collection $\mathfrak{h}$ of $k$- tuples $\mathcal{H}=(h_1,...,h_k)$ in such a way that the problem arising in Remark 6.6 can be avoided. Consider the first $k$ primes $p_1,...,p_k$ and for each $p_i$, let $e_i$ be the smallest natural number for which $1/p_i^{e_i}<d$. We will say $\mathcal{H} \in \mathfrak{h}$ if for each $p_i$, one has that $p_i^{e_i}|h_j \ \forall j \neq i$ and $(p_i, h_i)=1$.\\
In order to demonstrate a bound of the type stated in condition (\ref{averror1}), we will consider sums of the form
$$\sum_{\substack{n \leq N, n \equiv a \bmod q\\n+h_i \in \mathcal{A}_i }} 1,$$
with each $\mathcal{A}_i \subset \mathcal{A}$. The sets $\mathcal{A}_i$ will vary for each $N$ but we will ensure that $|\mathcal{A}_i|\gg_k N$ independently of $N$. Under these circumstances all of the arguments leading to the proofs of our main theorems go through.\\
Finally, we observe that in the case of a linear Bohr set $\mathcal{A}(N,d):=\left\{ n \leq N \ | \ \left\{n \alpha \right\} \in [0,d] \ \right\}$ we have the bonus of an additive structure, since $\mathcal{A}(N,d/2)+\mathcal{A}(N,d/2) \subset \mathcal{A}(N,d)$.       
\begin{proposition}\label{averageRa} Subject to the constraints described above, one has an inequality of the form
$$\sum_{q \leq Q} \left| \mathcal{R}_{\mathcal{A}} (N,a,q) \right| \ll N^{1-c(\rho)}$$
for $Q=N^{1/2}$ and some constant $c>0$ depending on $\rho$. 
\begin{proof}
Let $\varepsilon>0$ be a small but fixed constant to be chosen later. First we observe that the sequence $(\mathfrak{g}(n))_{n \in \N}$ is uniformly equidistributed $\bmod \ 1$. As a result, there exists an interval $I \subset [d/2,d]$ of length $|I|\leq d/3 N^{\varepsilon}$ containing no more than $N^{1-\varepsilon}$ elements of the form $\mathfrak{g}(n)$ with $n \leq N  $. We may now replace $d$ with the right end-point of $I$, since reducing $d$ by a factor $2$ will have no bearing on the statement of Theorem 13. Moreover, for each individual $i\leq k$ we shorten the interval $I$, obtained above, at its left-end and modify the constant term in $\mathfrak{g}$, so that we may assume $d=p_i^{- e_i}$. In other words, for each $i\leq k$ we will work on a shortened interval $I_i:=[0,d_i]$. Letting $\mathcal{A}_i$ denote the Bohr set associated to $I_i$, we note that 
$$\prod_{i=1}^k 1_{\mathcal{A}}(n+h_i)=\prod_{i=1}^k 1_{ [0,d_i]}(\mathfrak{g}(n+h_i))= \prod_{i=1}^k  \psi(\mathfrak{g}(n+h_i)) + O \left(\sum_{r=1}^k 1_{[d_i - \delta, d_i]}(\mathfrak{g}(n+h_i)) \right) ,$$

where $\delta=N^{- \varepsilon}$ and $\psi$ is the truncated Fourier series, as described in Remark\ \ref{fourierremark}. To deal with the remainder term we invoke the Erd\"os-Tur\'an-Koksma inequality (setting $H=N^{\varepsilon}$) and gather that 
\begin{align}\label{averesbohr}
\sum_{q \leq Q} \left| \mathcal{R}_{\mathcal{A}} (N,a,q) \right| & \ll \sum_{q \leq Q} \left|\sum_{\substack{|{\bf m}| \leq N^{\varepsilon}\\ {\bf m} \neq 0} } c_{\bf m } \sum_{n \leq N/q} e \left( \sum_{i=1}^k  \mathfrak{g}(nq + a +h_i) m_{i} \right)\right| \\
&+ \sum_{q \leq Q} \sum_{r=1}^{N^{\varepsilon}}    \frac{1}{r}	\left| \sum_{n \leq N/q} e \left(  \mathfrak{g}(nq + a +h_i) r \right)\right|  + O(N^{1-\varepsilon}) \notag ,
\end{align}
where the subscript ${\bf m}$ runs over $k$-tuples $(m_{1},..., m_{k}) $, the Fourier coefficients obey the bound $c_{{\bf m}}=\prod_{m_{i} \neq 0} c_{m_{i}}\ll \prod_{m_{i} \neq 0} |1/m_i|$ and $|{\bf m}|$ denotes the maximum norm of ${\bf m}$. Moreover, when $p_i^{e_i}| m_i $ for all component $m_i$ of ${\bf m}$, we may assume that $\sum_i m_i \neq 0 $. To see this, observe that when $p_i^{e_i} | m_i $, the Fourier coefficient $c_{{\bf m}} \leq 1/N^{k \varepsilon}$ and the number of $k$-tuples ${\bf m}$ for which $\sum_i m_i = 0 $ is at most $N^{(k-1) \varepsilon}$.\\

Due to our choice of $\mathcal{H}$, we now see that the polynomial $P(x):=\sum_{i=1}^k  \mathfrak{g}(xq + a +h_i) m_{i}$ appearing in (\ref{averesbohr}) does not collapse to a constant function and in fact has a leading coefficient of the form $\beta(q, {\bf m}):={{D}\choose{j}} q^{D-j} \alpha_{k-j} (\sum_{i=1}^k h_{i}^j m_i) $ for $j=0 \text{ or } 1$. We will deal with the first sum in (\ref{averesbohr}) and note that the second quantity is treated in very much the same way.
As a consequence of Weyl's inequality, we have 

$$\left| \sum_{n \leq N/q} e(P(n))\right|\ll \left( \frac{N}{q} \right)^{1-c(\varepsilon)}$$
as long as one can find natural numbers $N^{\varepsilon} \leq \nu \leq N^{D-\varepsilon}$ and $r\geq 1$ coprime to $\nu$, satisfying 
$$ \left|\beta(q, {\bf m}) -\frac{r}{\nu} \right|\leq \frac{1}{\nu^2}.$$  

First consider those ``good" pairs $(q,{\bf m})$ for which $\beta(q, {\bf m})$ admits a rational approximation as above. Applying Weyl's inequality to these pairs, it follows that their contribution to the RHS of (\ref{averesbohr}) is easily subsumed in the desired bound. Also, since $\alpha$ is Diophantine, we can choose $\varepsilon$ in such a way that $(q,{\bf m})$ is good whenever $q^D  \leq N^{\varepsilon}$.\\
Let us arrange the remaining ``bad" pairs $(q,{\bf m})$ into dyadic cubes $\mathcal{Q}\leq q \leq 2\mathcal{Q}$ and $H\leq |{\bf m}| \leq 2H$. By Dirichlet's theorem, for each such bad pair, we thus find a $\nu \leq N^{c(\varepsilon)/2}$ such that $\left|\beta -r/ \nu \right|\leq 1/\nu N^{D-\varepsilon}$. In other words, we can write
$$q \in \bigcup_{\nu \leq N^{c(\varepsilon)/2}} \left\{\mathfrak{q} \leq \mathcal{Q}| \ ||\mathfrak{q}^D  |{\bf m}| \alpha ||\leq \frac{1}{N^{D-\varepsilon}}  \right\}.$$

From the previous lemma we gather that each set in the union above has size $O \left(Q^{1-c(\rho)} \right)$ (for $\varepsilon$ sufficiently small) so that
$$\mathop{  \sum_{\mathcal{Q} \leq q \leq 2\mathcal{Q} } \sum_{H \leq |{\bf m}| \leq 2 H }  }_{(q,{\bf m}) \ \text{bad}} c_{{\bf m}} \left| \sum_{n \leq N/q} e(P(n))\right| \ll (\log N) Q^{1-c(\rho)} \frac{N}{\mathcal{Q}}. $$
The proof is completed by summing over all dyadic cubes for which $\mathcal{Q} \geq N^{\varepsilon/D}$.\\
When $D=1$, the proof proceeds in the same manner except it is no longer necessary to sum over $k$-tuples ${\bf m}$ in equation (\ref{averesbohr}).  
\end{proof}
\end{proposition}

\begin{lemma}\label{smallmodBohr} For $Q=(\log N)^{A}$, we have
$$\sum_{q \leq Q} \frac{1}{\varphi(q)} \sum_{\chi \bmod q}^{*}  \left| \psi_{\mathcal{A}}' (N,\chi) \right| \ll \frac{N}{(\log N)^{A}}$$

\begin{proof}
Let $\varepsilon>0$ be fixed and expand $\psi_{\mathcal{A}}' (N,\chi)$ using the smooth cut-off function discussed in the proof of Proposition \ref{averageRa}. For $\chi\neq \chi_0$ this yields

\begin{align*}
\psi_{\mathcal{A}}' (N,\chi)&=\sum_{|{\bf m}| \leq N^{\varepsilon}} c_{\bf m} \sum_{n \leq N} \chi(n) \Lambda(n) e \left(\sum_{i=1}^k  \mathfrak{g}(n+h_i) m_{i} \right) +O(N^{1- \varepsilon})\\
& =\sum_{|{\bf m}| \leq N^{\varepsilon}}  \frac{c_{\bf m}}{\tau(\overline{\chi})}  \sum_{a=1}^{q} \sum_{n \leq N} \overline{\chi(a)}  \Lambda(n) e \left( \sum_{i=1}^k  \mathfrak{g}(n +h_i) m_{i} + \frac{an}{q} \right) +O(N^{1- \varepsilon})
\end{align*}
where $\tau(\chi)=\sum_{m=1}^{q} \chi(m)e(m/q)$, $c_{\bf m}\ll 1/|{\bf m}|$ and $H \ll N^{\varepsilon}$. First consider the case $D=1$. Since $\alpha$ is Diophantine, we may estimate the exponential sum immediately, using Theorem \ref{Harman} to obtain the desired bound.\\
When $D\geq 2$, we adopt the exact same strategy as in the proof of Proposition \ref{averageRa}. Organise ${\bf m}$ into dyadic cubes and then consider separately the good and bad pairs $(q,{\bf m})$.\\
Finally, when ${\bf m}=0$, we have the classical estimate
$$\sum_{n \leq N} \chi(n) \Lambda(n) \ll \frac{N}{(\log N)^{2A}}.$$ 

The case $\chi=\chi_0$ is treated similarly.
\end{proof}
\end{lemma}

Our final task is to prove the bilinear form estimate (\ref{bilinear1}). We recall the following variant of the classical large sieve inequality
\begin{align}\label{claslargesieve}
\sum_{q \leq Q} \frac{q}{\varphi(q)} \sum_{\chi}^{*}  \max _{u} \left|\mathop{\sum_{m \leq M} \sum_{l\leq L} }_{mn \leq u} a_m b_l \chi(ml) \right| \ll (M+Q^2)^{1/2}(L+Q^2)^{1/2}(\sum_{m \leq M} |a_m|^2)^{1/2}(\sum_{l \leq L} |b_l|^2)^{1/2} \log(2ML). 
\end{align}

\begin{lemma}
\begin{align}\label{bilinearbohr} 
\sum_{q \leq Q} \frac{q}{\varphi(q) }\sum_{\chi}^{*} \left|\mathop{\sum_{m \leq M} \sum_{l\leq L} }_{ml + \mathcal{H} \subset \mathcal{A}} a_m b_l \chi(ml)  \right| \ll  (M+Q^2)^{1/2}(L+Q^2)^{1/2}(\log ML)^3 \left\|a\right\|_2 \left\|b\right\|_2. 
\end{align}

\begin{proof}
We will assume without loss of generality that $L\geq M$. The inner-most double sum in (\ref{bilinearbohr}) can be expanded to obtain 
\begin{align}\label{SumS}
\mathcal{S}&:=\mathop{\sum_{m \leq M} \sum_{l\leq L} }_{ml + \mathcal{H} \subset \mathcal{A}} a_m b_l \chi(ml) =\sum_{m \leq M} \sum_{l\leq L} a_m b_l \chi(ml) \sum_{|{\bf r}| \leq L^{C_1} } c_{\bf r }  \ e\left( \sum_{i=1}^k  \mathfrak{g}(ml +h_i) r_{i} \right)+O(1/L^2) \notag \\
&= \sum_{|{\bf r}| \leq L^{C_1} } c_{\bf r } \sum_{m \leq M} \sum_{l\leq L} a_m b_l \chi(ml) \left[\int_{(ML)^2}^{F_{\bf r}(ml)}  \ e\left( x \right) \ dx +e\left((ML)^2\right) \right]+O(1/L^2).
\end{align}  
Here $F_{\bf r}(x)=\sum_{i=1}^k  \mathfrak{g}(x +h_i) r_{i}$ and we shift the leading coefficient by a large integer multiple of $2 \pi i$ to ensure that $F_{\bf r}(y) \geq (ML)^2$ for $y\geq 1$. Furthermore, the constant $C_1$ depends only on $\eta$, ${\bf r }=(r_1,...,r_K) \in \Z^K$ and regarding the Fourier coefficients, we have $c_{\bf r } \ll \prod_{r_{i} \neq 0} |1/r_i|$. Interchanging the order of summation and integration in (\ref{SumS}) leads to the quantity
\begin{align}\label{SumSexp}
\int_{(ML)^2}^{L^{C_2(\eta,K)}} e(x) \mathop{\sum_{m \leq M} \sum_{l\leq L} }_{F_{\bf r}(ml) > x} a_m b_l \chi(ml)
\end{align}
The condition $F_{\bf r}(ml) > x$ can be replaced by an expression of the form $ml \in \cup_i [Z_i(F_{\bf r}-x), Z_{i+1}(F_{\bf r}-x)]$ where $Z_i(G)$ denotes the $i$-th real zero of the function $G$. It therefore suffices to evaluate (\ref{SumSexp}) under the new condition $Z_i(F_{\bf r}-x) < ml$ or rather $H(x):=\log(Z_i(F_{\bf r}-x)) < \log(ml)$. A simple application of the Implicit Function Theorem shows that $H'(x)\ll 1/x$, with the implied constant being independent of the modifications made to $F_{\bf r}$. Substituting the discontinuous integral (p. 165 of \cite{Dav}) 
\begin{equation*}
\label{eq:restr-def}
2\int_{-T}^{T} e(\alpha t) \sin(2 \pi \beta t) \frac{dt}{t}=\left\{
\begin{array}{ll}
1 + O(T^{-1} (\beta - |\alpha|)^{-1}) & \mathrm{if }\  |\alpha| < \beta\\
O(T^{-1} (\beta - |\alpha|)^{-1}) & \mathrm{if}\ |\alpha| > \beta\\
\end{array}
\right.
\end{equation*}
into $(\ref{SumS})$, we are left with sums of the form
\begin{align*}
\int_{-T}^{T} \sin(2 \pi [\log(m) + \log(l)] t)  \int_{(ML)^2}^{L^{C_2(\eta,K)}} e(x+t H(x) )  \sum_{m \leq M} \sum_{l\leq L} a_m b_l \chi(ml) \ dx \ \frac{dt}{t}.
\end{align*}
Choosing $T=(ML)^{3/2}$, we see that $\frac{d}{dx} (x+t H(x))$ is bounded away from zero and piecewise monotonic so that van der Corput's lemma can be applied. Note that, after applying a simple trig identity, we have effectively separated the variables $m$ and $l$ so that the resulting quantities may be absorbed into $a_m$ and $b_l$. We may now proceed to use (\ref{claslargesieve}).
\end{proof}
\end{lemma}

%\subsubsection{Some closing remarks about linear Bohr sets}

{\bf Acknowledgements} I would like to thank my advisor, Prof. Terence Tao, for his generous encouragement and numerous suggestions throughout this project.

\end{document}